\documentclass[11pt,reqno,a4paper]{amsart}
\usepackage[a4paper,vmargin={3.5cm,3.5cm},hmargin={2.5cm,2.5cm},bottom=30mm]{geometry}
\usepackage{amsfonts}
\usepackage{amsthm,amsmath,slashbox}
\usepackage{amsbsy}
\usepackage{bm}
\usepackage{pxfonts}
\usepackage{dsfont}
\usepackage{amssymb} 
\usepackage[dvipsnames]{xcolor}

\usepackage[mathscr]{eucal}

\renewcommand{\imath}{{\mathbbm{i}}}

\usepackage{units}
\usepackage{enumerate,enumitem}
\usepackage{bbm}
\usepackage{exscale}
  \numberwithin{equation}{section} 
  \date{}

\theoremstyle{plain} 
    \newtheorem{theorem}{Theorem}
    \newtheorem{lemma}[theorem]{Lemma}
    \newtheorem{proposition}[theorem]{Proposition}
   
    \newtheorem{step}{Step}
    \newtheorem{corollary}[theorem]{Corollary}
    \newtheorem{claim}[theorem]{Claim}

\theoremstyle{definition} 
    \newtheorem{definition}{Definition}

    \newtheorem{remark}[definition]{Remark}
    
    \newtheorem{assumption}{Assumption}
\usepackage{appendix}
\renewcommand\appendix{\par
        \renewcommand\thesection{A}
        \renewcommand\thesubsection{A\arabic{subsection}}
        \renewcommand\thetable{A\arabic{table}}}

\DeclareMathOperator{\T}{\textbf{T}}
\DeclareMathOperator{\R}{\mathbb{R}}

\DeclareMathOperator{\Z}{\mathbb{Z}}
\DeclareMathOperator{\N}{\mathbb{N}}
\usepackage{mathtools}

\DeclareMathOperator{\De}{d}
\DeclareMathOperator{\one}{\mathbbm{1}} 
\newcommand{\E}{\mathsf{E}}

\newcommand{\prob}{\mathsf{P}}
\renewcommand{\T}{\mathbb{T}}

\DeclareMathOperator{\e}{e}
\newcommand{\la}{\left\langle}
\newcommand{\ra}{\right\rangle}

\newcommand{\eps}{\epsilon}
\newcommand{\eq}[1]{\begin{equation#1}}
\newcommand{\eeq}[1]{\end{equation#1}}

\usepackage[comma, sort&compress]{natbib}
\usepackage{hyperref}
\hypersetup{ colorlinks=true, linkcolor=blue, citecolor=blue,
  filecolor=blue, urlcolor=blue}
\begin{document}

\title[Divisible sandpile with heavy-tailed variables]{\bf \textsc{The divisible sandpile with heavy-tailed variables}}
\author[A. Cipriani, R. S. Hazra and W. M. Ruszel]{Alessandra Cipriani, Rajat Subhra Hazra and Wioletta M. Ruszel}
\address{WIAS Berlin.}
\email{Alessandra.Cipriani@wias-berlin.de}

\address{Indian Statistical Institute, Kolkata.}
\email{rajatmaths@gmail.com}
\address{TU Delft.}
\email{W.M.Ruszel@tudelft.nl}


\begin{abstract}
This work deals with the divisible sandpile model when an initial configuration sampled from a heavy-tailed distribution. Extending results of \cite{LMPU} and \cite{CHR2016} we determine sufficient conditions for stabilization and non-stabilization on infinite graphs. We determine furthermore that the scaling limit of the odometer on the torus is an $\alpha$-stable random distribution.
\end{abstract}
\maketitle
\section{Introduction}
The divisible sandpile model, a continuous version of the (discrete) abelian sandpile model (ASM) was introduced by \cite{LePe10, LevPer} to study scaling limits of the rotor aggregation and internal DLA growth models. 

The basic mechanism in these models is that to each site of some graph there is associated a {\em height} or {\em mass}. If the height exceeds a certain value then it collapses by distributing the excess mass (uniformly) to the neighbours which can then result in a series of cascades. \\
One of the questions arising for these cascading models is the dichotomy between stabilizing and exploding configurations. 

For the ASM \cite{FMR2009} showed that given an initial i.i.d.~configuration on $\Z^d$ the model will stabilize almost surely, depending solely on the mean density at a fixed site and the dimension $d$. In \cite{LMPU} the authors extended this study to the divisible sandpile model on general vertex-transitive graphs. One of their results deals with the characterization of this dichotomy according to the mean height and transience resp. recurrence of the graph (and not anymore on $d$). If the mean height is larger than 1 then almost surely the initial configuration does not stabilize whilst a value smaller than 1 ensures stabilizability. At the critical value 1 under the additional assumption of finite variance the model does not stabilize. 

The proof of non-stabilizability at the critical value involves studying a so-called \textit{odometer function}. It measures the amount of mass emitted from a site during stabilization. 
\cite{LMPU} study the expected odometer growth in the case of an initial Gaussian configuration using an interesting connection with the discrete bilaplacian Gaussian field. The discrete bilaplacian (or membrane) model is a particular random interface model (similar to the Gaussian Free Field) and was introduced in the mathematics literature by \cite{Saka}, \cite{Kurt_d4,Kurt_d5}. Levine and coauthors conjectured that the rescaled odometer converges to a continuum bilaplacian field when the mesh size of the discrete torus becomes finer.  

In \cite{CHR2016}  the authors considered a general divisible sandpile model with i.i.d. initial distribution on a discrete torus and proved the conjecture of \cite{LMPU} on the torus $\T^d$ determining the limiting field. 

 In this article we are interested in exploring the properties of the divisible sandpile model when the initial mass comes from heavy-tailed distributions.
We are interested in extending results from both \cite{CHR2016} and \cite{LMPU}, namely we first study the dichotomy between stabilizing versus exploding configurations and secondly determine the scaling limit of the odometer function for heavy-tailed distributions on the torus. The novelty of the article is to consider the stabilization versus explosion dichotomy for divisible sandpiles for more general initial distributions by removing the finite variance assumption at the critical value $\E(s)=1$ and to study scaling limits for those generalized random variables. To the authors' knowledge this is the first result constructing an $\alpha$-stable random distribution on the torus.\newline

More precisely, the divisible sandpile of a locally finite, undirected, connected graph $G=(V, E)$ is defined as follows: start with an initial configuration $s: V\to \R$. A vertex $x$ is unstable if its height $s(x)>1$ and stable otherwise. At the first time instance all unstable vertices $x$ topple keeping mass $1$ to themselves and redistributing the excess $s(x)-1$ equally among their neighbours. 
If at time $n$ the total mass distributed from $x$ is given by $u_n(x)$, then it can be proved that $u_n\to u$ where $u: V\to [0,\,+\infty]$. $u$ is called the {\em odometer} for the configuration $s$; if the odometer is finite for all $x\in V$ then we say that a configuration is stable. In \cite{LMPU} many properties of the divisible sandpile were studied when $(s(x))_{x\in V}$ are independent and identically distributed random variables with finite mean and finite variance. It then becomes a natural question to see if their analysis can be pushed further to more general random variables, especially when mean and variance are infinite. In particular, we shall see that the finiteness of the mean is not necessary to study the dichotomy of stabilization versus explosion.

We will consider initial heights which are regularly varying with index $\alpha$, i.e. they satisfy
\begin{equation}\label{eq:intro:rv}
\prob(|s|> t) \sim t^{-\alpha} L(t)\quad \text{ as } t\to + \infty
\end{equation}
where $L$ is some slowly varying function and $\alpha\in (0,2]$. Such variables arise naturally when one considers domain of attractions of stable distributions. 

We show that the initial configuration almost surely will not stabilize if $\E(s) \in (1,\infty ]$  or if $\E(s)=1$, assuming infinite variance and some additional property of the underlying graph and $\alpha$. On the other hand the initial configuration will  stabilize almost surely if $\E(s) \in [-\infty, 1)$.
It is tempting to consider the value of $\alpha$  in~\eqref{eq:intro:rv} as a parameter which is in some sense tuning the dichotomy, since it is related to finiteness resp. infiniteness of the first and second moment. If $\alpha \in (0,1)$ then the mean $\E(s)=\pm \infty$ whereas for $\alpha \in (1,2)$ we know that $\E(s)<\infty$ and the variance is infinite. However in the boundary cases $\alpha=1$ and $\alpha=2$ the finiteness of the moments depends on the function $L$, hence we cannot decide a priori whether the configuration is stabilizable or not knowing solely $\alpha$. 

A second part of this paper focuses on a special finite connected graph, the discrete torus. In general on a finite graph $G$ with $|V|=n$ and for which the mass is conserved, that is, $\sum_{x\in V} s(x)=n$, the system stabilizes to the configuration constantly equal to $1$.  This regime corresponds to the critical case when $\E(s)=1$. The odometer $u$ satisfies the following discrete equation \cite[Lemma~7.1]{LMPU}:
\begin{equation}\label{eq:odometer}
\begin{cases}
\Delta u(x)&= 1-s(x) \\
 \min_{x\in V} u(x)&=0
\end{cases},
\end{equation}
where $\Delta$ is the discrete Laplacian. When $V$ is the discrete torus of side length $n$ (denoted by $\Z_n^d$) the study of the scaling limit of the odometer becomes interesting. 
We construct a new field on the dual of $C^\infty(\T^d)$ (the space of smooth functions on the torus) to which we show that the rescaled odometer converges. This field belongs to the class of $\alpha$-stable generalised random fields, which is a natural extension of Gaussian random fields. It is remarkable that the sandpile is able to span through a whole class of generalised fields which all have the stability property like stable random variables. \newline

\paragraph*{{\it Outline of the article}} The article is structured as follows: in Section~\ref{sec:setup_main} we give the basic definitions and explain rigorously the results obtained. In Section~\ref{sec:stabil} we deal with the proofs of the results concerning stabilization on infinite graphs. In Section~\ref{sec:odometer_limit} we determine the scaling limit of the odometer on the discrete torus. Auxiliary results are proved in Appendix \ref{appendix}.
\newline
\paragraph*{{\it Acknowledgements}}We are grateful to Mark Veraar for helpful discussions. The second author also would like to thank Deepak Dhar for an enlightening discussion on sandpile models. The first author's research was partially supported by the Dutch stochastics cluster STAR (Stochastics -- Theoretical and Applied Research). The second author's research was supported by Cumulative Professional Development Allowance from Ministry of Human Resource Development, Government of India and Department of Science and Technology, Inspire funds.
\section{Basic setup and main results}\label{sec:setup_main}
\subsection{Notation}We start with some preliminary notations which are needed throughout the paper. 
Let $\T^d$ be the $d$-dimensional torus, viewed as ${\R^d}/{\Z^d}$ or as $[-1/2,\,1/2)^d\subset\R^d$ alternatively. The discrete torus of side-length $n\in \N$ is $\Z_n^d:=[-{n}/{2},\,{n}/{2}]^d\cap \Z^d$, and $\T_n^d:=[-1/2,\,1/2]^d\cap (n^{-1}\Z)^d$ is the discretization of $\T^d$. {For a discrete set $V$ we denote as $|V|$ its cardinality.} Moreover let $B(z,\,\rho)$ be a ball centered at $z$ of radius $\rho>0$ in the $\ell_\infty$-metric. We will use throughout the notation $z\cdot w$ for the Euclidean scalar product between $z,\,w\in \R^d$. With $\|\cdot\|_\infty$ we mean the $\ell_\infty$-norm, and with $\|\cdot\|$ the Euclidean norm. We will let $C,\,c$ be positive constants which may change from line to line within the same equation. We define the Fourier transform of a function $u\in L^1(\T^d)$ as $\widehat u(y):=\int_{\T^d}u(z)\exp\left(-2\pi\imath y\cdot z\right)\De z$ for $y\in \Z^d$. We will use the symbol $\,\widehat \cdot\,$ to denote also Fourier transforms on $\Z_n^d$ and $\R^d$. We say a function has mean zero if $\int_{\T^d} f(z)\De z=0$. We will denote for a real-valued random variable $X$ and $x\in \R$
\begin{equation}
F_X(x):=\prob(X\le x),\quad \overline F_X(x):=1-F_X(x)=\prob(X>x).
\end{equation}
We write for two positive functions $f,\,g$
\[
f(x)\sim g(x)\quad\text{as }x\to x_0
\]
if $\lim_{x\to x_0}f(x)/g(x)=1.$
\subsection{Assumptions on the configuration}
We recall here the definition of regularly varying function: a non-negative random variable $X$ is called {\em regularly varying} of index $\alpha\ge 0$, and we write $X\in RV_{-\alpha}$, if
\[
\overline F_X(x)\sim x^{-\alpha} L(x) \quad\text{as }x\to +\infty 
\]
where $L$ is a slowly varying function, i. e.,
\[
\lim_{x\to+\infty}\frac{L(tx)}{L(x)}=1\quad\text{for all }t>0.
\]
We recall the definition of variables in the domain of attraction of a stable distribution:
\begin{definition}[Domain of normal attraction of stable variables]
\label{firstassumption}
Let $\alpha\in (0,\,2]$. Let $V$ be a countably infinite index set and $(W(x))_{x\in V}$ be i.i.d. symmetric random variables with common distribution function in the domain of normal attraction of an $\alpha$-stable distribution. This means that, for $V_1\subset V_2\subset\ldots$ such that $\cup_{k\ge 1}V_k=V$, we have the following limit:
\begin{equation}\label{eq:doa:stable}
\lim_{k\to+\infty}{|V_k|^{-\frac{1}{\alpha}}}\sum_{x\in V_k} W(x)\stackrel{d}{=} \rho_\alpha,
\end{equation}
where $\rho_\alpha$ has a symmetric $\alpha$-stable law which we denote as $S\alpha S(\mathfrak c)$, that is, $\E[ \exp(\imath \theta \rho_\alpha)]= \exp({-\mathfrak{c}^\alpha|\theta|^\alpha})$ for some $\theta\in \R$. 
\end{definition}
In our work we will often use this definition setting $V:=\Z_n^d$ (it will be clear from the context when). If the scale parameter of the $\alpha$-stable law is $1$, we will write $\sigma(x)\stackrel{d}{=} S\alpha S(1).$ If this happens, it is well known that $|\sigma(x)|$ has a regularly varying tail with index $-\alpha$, for $\alpha\in (0,2]$.
\begin{remark}
The results we are going to prove can be extended to a more general set-up assuming further necessary and sufficient conditions for the $(\sigma(x))_{x\in V}$ to be in the domain of attraction of stable variables (classical references on the topic are \citet{MikReg,ST}). However to keep the exposition accessible without harming the mathematical aspects we assume the simpler Definition~\ref{firstassumption}.
\end{remark}
\subsection{Stability on infinite graphs: beyond finite variance}
First we shall see some properties of divisible sandpiles on infinite graphs. More specifically we consider $G= (V, E)$ to be an infinite vertex transitive graph. Let $\Gamma\subset \mathrm{Aut}(G)$ be a subgroup which acts transitively on $V$ and let $\prob$ be a $\Gamma$-invariant probability measure. Let $o$ be a distinguished vertex of $V$ which we keep fixed. Denote by {$\R^V$} the set of divisible sandpile configurations on G. Recall from \citet[Section 2]{LMPU} that in toppling procedure starting from an initial configuration {$s\in \R^V$}, the total mass emitted by a site $x\in V$ to each of its neighbours during the time interval $[0,\, n]$ is $u_n(x)$, so that the resulting configuration at time $n$ is $s_n=s+\Delta u_n$. In the same work it is shown that if $u$ is a finite toppling procedure then $s_\infty= \lim_{t\to \sup{T}} s_t$ exists, where $T$ is a well-ordered set of toppling times. A toppling procedure $u$ is called stabilizing for $s$ if $u$ is finite and $s_\infty(x)\le 1$ for all $x\in V$ . One says that $s$ stabilizes if there exists a stabilizing toppling procedure for $s$.

Our first Theorem tries to explore the case when initial configurations does not necessarily have finite mean. 

\begin{lemma}\label{lemma:infinite mean}
Let $G=(V, E)$ and $\prob$ be as above. Let $(s(x))_{x\in V}$ be i.i.d.
\begin{enumerate}[label=(\roman*)]
\item\label{item:(i)} If $\E[s(o)]=+\infty$, then $\prob( s\text{ stabilizes})=0$. 
\item\label{item:(ii)} {If $\E[s(o)]=-\infty$, then $\prob( s\text{ stabilizes})=1.$}
\end{enumerate}
\end{lemma} 
Recall that if $X$ is a (non-negative) regularly varying random variable with index $-\alpha$ then $\E\left[ X^\beta\right]<+\infty$ for $\beta<\alpha$ and $E\left[X^\beta\right]=+\infty$ for $\beta>\alpha$. At $\beta=\alpha$ the mean may be finite or infinite.
Note that when one assumes that $s$ has a regularly varying tail of index $-\alpha$ with $\alpha<1$ then the above result implies that there is no stability almost surely, since $s$ has infinite mean. 

The configurations for which the mean is finite (but not necessarily the variance) require some more analysis. Again note that if $s$ has a regularly varying tail of index $-\alpha$ with $\alpha\in (1,2)$, the mean is finite. Recall also that the cases $\E[s(o)]<1$ and $\E[s(o)]\in (1,\,+\infty)$ can be dealt with the results from \citet[Lemmas~4.1, 4.2]{LMPU}. When the mean is $1$ we must study the dependence on the underlying graph more closely, and in particular the behavior of the simple random walk on it. \cite{LMPU} show that there is no stability adopting different techniques according to the transience or recurrence of the graph, and we will adopt a somewhat similar viewpoint for regularly varying variables. Let us recall the Green's function $g(x,\,y):= \sum_{j=0}^{+\infty} P_x( S_j=y)$, $x,\,y\in V$, where $S_j$ is the simple random walk on $V$ started at $x$. We split the critical case into two broad cases: $\sum_{x\in V} g(o, x)^\alpha=+\infty$ and $\sum_{x\in V} g(0, x)^\alpha<+\infty$. With a bit of abuse of nomenclature we call the first case $\alpha$-singly transient and the second case $\alpha$-doubly transient (for a summary of stabilizability, see Table~\ref{table}). In the following results, since $s(x)$ has mean $1$, we impose conditions on the recentered variable $s(x)-1$, as it is natural to assume symmetry. Given that $s(x)-1$ and $s(x)$ are tail equivalent in the case of regular variation, this does not effect the outcome of the result.

\begin{theorem}[{$\alpha$-singly transient}]\label{thm:singly_trans}
Let $(s(x))_{x\in V}$ be a divisible sandpile on an infinite vertex transitive graph $G=(V, E)$ such that $(Y(x))_{x\in V}:=(s(x)-1)_{x\in V}$ are i.i.d., zero-mean, symmetric random variables in the normal domain of attraction of a $S\alpha S$ random variable with $\alpha\in [1,2)$ (recall Definition~\ref{firstassumption}). Suppose $g(o,\,y)<+\infty$ uniformly for all $y\in V$ and
\begin{equation}\label{eq:assumption_g1}\sum_{y\in V} g(0,y)^\alpha =+\infty.
\end{equation}
Then $\prob(s\text{ stabilizes})=0$. 
\end{theorem}

\begin{remark}\label{rem:sing_trans_alpha}
In the case in which $V:=\Z^d$, $d\ge 3$, then by \citet[Theorem~4.3.1]{LawlerLimic} we obtain that \eqref{eq:assumption_g1} is satisfied if $\alpha\le{d}/({d-2})$. In particular this implies that the singly transient case for the square lattice corresponds to $\alpha\in (1,\,{d}/({d-2}))$, hence it comprises the cases $d=3,\,4$.
\end{remark}

Now we deal with the case $\alpha$-doubly transient case. Although in this case one may expect to assume $\sum_{y\in V} g(o, \, y)^\alpha<+\infty$, we shall assume something stronger to prove our results. 
\begin{assumption}\label{thirdassumption}
Assume that
\begin{enumerate}[label=({\alph*})]
\item $(Y(x))_{x\in V}:=(s(x)-1)_{x\in V}$ are i.i.d., zero-mean, symmetric random variables in the normal domain of attraction of a $S \alpha S (1)$ random variable with $\alpha\in (1,2]$.
\item \label{item:(b)}There exists $\delta\in (1,\,\alpha)$ such that
 \eq{*}
\sum_{y\in V}g(o,\,y)^\delta<+\infty.
\eeq{*}
\end{enumerate}
\end{assumption}
Then {we can state the following}
\begin{theorem}[$\alpha$-doubly transient case]\label{thm:doubly_trans}
Let $G=(V,\,E)$ be an infinite vertex transitive graph and let $(s(x))_{x\in V}$ such that they satisfy Assumption~\ref{thirdassumption}. Then $\prob( s\,  \text{ stabilizes})=0$.
\end{theorem}

We note that~\ref{item:(b)} implies that $\sum_{y\in V} g(0,y)^\alpha<+\infty$. In fact, we will deal with infinite series of the form $\sum_{x\in V} g(0,x) Y(x)$ which converge almost surely when one assumes~\ref{item:(b)}. Such assumptions are well-known in heavy-tailed time series literature.  The series also converges if one assumes additional conditions on slowly varying functions (see \citet[Lemma~A.4]{mikosch:gena:2000} for these conditions). For example, if $Y(o)$ satisfies $F_Y(t)=t^{-\alpha}$, then one can relax the assumption \ref{item:(b)} and choose $\delta:=\alpha$ to obtain the statement of Theorem~\ref{thm:doubly_trans}. 

\begin{remark}
Analogously to Remark~\ref{rem:sing_trans_alpha}, one can show that for the graph $\Z^d$, $d\ge 5$, an exponent $\delta<\alpha$ such that \ref{item:(b)} holds can always be found (indeed one needs $d/(d-2)<\delta<\alpha)$).
\end{remark}

\begin{table}[h!]
\centering
\resizebox{\textwidth}{!}{\begin{tabular}{ |c|c|c|c|c|c|c|  }
 \hline
  \backslashbox{$\mathsf{Var}[s(o)]$}{$\E[s(o)]$}&$[-\infty,\,1)$   & $1$  &$ (1,\,+\infty)$ & $+\infty$&$-\infty$ \\
 
  \hline
  Finite &$1$& $0$ &$0$&{$\times$}&$\times$\\ 
  \hline
 Infinite &$1$ : Lemma~\ref{lemma:infinite mean} &$0:\,\begin{cases} \alpha-\text{singly transient (Thm.~\ref{thm:singly_trans})}\\ \alpha-\text{doubly transient (Thm.~\ref{thm:doubly_trans})} \end{cases}$ &$0$& $0$: Lemma~\ref{lemma:infinite mean}& 1: Lemma~\ref{lemma:infinite mean} \\
 \hline
\end{tabular}}\caption{Summary of stabilizability. In each cell we write the value of $\prob( s\,  \text{ stabilizes})$.}\label{table}
\end{table}

This completes the picture of stability on a divisible sandpile for regularly varying random variables. We now explore the odometer behavior on the finite graphs, and specifically on a torus. 
\subsection{Scaling limit of the odometer on the torus}
For a finite connected graph, the divisible sandpile is stable if and only $\sum_{x\in V} s(x)\le |V|$. When the sum is exactly $|V|$ the configuration stabilizes to the all $1$ configuration and the odometer $u$ is the unique function $u$ which satisfies~\eqref{eq:odometer} (\citet[Lemma 7.1]{LMPU}). This equation can be useful in determining the representation of the odometer. One can obtain the following result, for which we do not give a proof since it mimicks closely that of \citet[Proposition~1.3]{LMPU}.
\begin{proposition}\label{prop:1.3LMPU}
Consider $G=(V, E)$ a finite connected graph. Let {$s(x)$ be a configuration such that $\sum_{x\in V}s(x)=|V|$.}
Then the configuration stabilizes to the all $1$ configuration and the distribution of the odometer $u$ is given by
$$( u(x))_{x\in V}= \left( v(x)- \min_{z\in V} v(z)\right)_{x\in V}$$
where
\begin{equation}
v(x)= \frac1{\mathrm{deg}(x)} \sum_{w\in V} g(w,x)(s(w)-1)
\end{equation}
and $g(w,x)={|V|}^{-1} \sum_{z\in V} g^z(w, x)$, {where $g^z(x,y)$ is the expected number of visits to $y$ by a simple random walk started at $x$ before hitting $z$.}
\end{proposition}

When $(\sigma(x))_{x\in V}$ are i.i.d. Gaussians and
\eq{}\label{eq:def_conf_graph}
s(x)=1+\sigma(x)-\frac{1}{|V|}\sum_{w\in V}\sigma(w),\quad x\in V
\eeq{}
then one can show that $v(x)$ is distributed as a discrete bilaplacian field on the torus, that is, it is a centered Gaussian field with covariance given by 
$$\E[v(x) v(y)]= \frac1{\mathrm{deg}(x)\mathrm{deg}(y)} \sum_{w\in V} g(x,w)g(w,y).$$
In this Gaussian case, this hints at the possibility that the field $u$, appropriately rescaled, may converge to the continuum bilaplacian field on the torus. To describe the general case, let us consider the interpolated rescaled odometer:
$$\Xi_n(x):= 4\pi^2 n^{d-\frac{d}{\alpha}-2}\sum_{z\in \T_n^d} u(nz) \one_{B\left(z, \,\frac{1}{2n}\right)} (x).$$
For $f\in C^\infty(\T^d)$ and mean zero we can define the action of the field $\Xi_n$ on $f$ as
$$\la \Xi_n, f \ra = {4\pi^2}n^{d-\frac{d}{\alpha}-2}  \sum_{z\in \T_n^d} u(nz) \int_{B\left(z,\, \frac{1}{2n}\right)} f(t) \De t.$$
\begin{theorem}\label{thm:scaling_RV}
Let $d\ge 1$. Let $(\sigma(x))_{x\in \Z_n^d}$ be i.i.d. and satisfy Definition \ref{firstassumption}, and furthermore let $(s(x))_{x\in \Z_n^d}$ as in \eqref{eq:def_conf_graph} where $V:=\Z_n^d$. There exists a random distribution $\Xi_\alpha$ on $(C^\infty(\T^d))^\ast$ such that: for all $m\in \N$ and $f_1, f_2,\ldots,\, f_m\in C^\infty(\T^d)$ with mean zero, the random variables $\la \Xi_n, f_j \ra$ converge jointly in distribution to a random variable $\la \Xi_\alpha, f_j\ra$. Moreover, the characteristic functional of $\Xi_\alpha$ is given by 
\begin{equation}\label{eq:limitfield}
\E[\exp(\imath \la \Xi_{\alpha}, f\ra )]=\exp\left(-\int_{\T^d} \left| \sum_{z\in \Z^d\setminus \{0\}} \frac{\exp(-2\pi \imath z\cdot x)}{\|z\|^2} \widehat{f}(z)\right|^{\alpha} \De x\right).
\end{equation}
\end{theorem}

The above theorem describes the finite dimensional convergence of the odometer field. The limiting characteristic function is well-defined and indeed defines an $\alpha$-stable cylindrical random field, of which we recall the definition. Let ``$\sim$'' be the equivalence relation that identifies two functions differing by a constant and call $\mathcal T:=C^\infty(\T^d)/{\sim}$. Let $\alpha\in (0,\,2].$ A random variable $\Xi_\alpha$ on $\mathcal T^*$ is called $\alpha$-stable if, given two independent copies $\Xi_{\alpha,\,1}$ and $\Xi_{\alpha,\,2}$ of $\Xi_\alpha$, then for any $a,\, b>0$ and $f \in \mathcal T$
\protect\eq{*}\E[\exp(\imath \la \Xi_{\alpha,\,1}, af\ra )]\E[\exp(\imath \la \Xi_{\alpha,\,2}, bf\ra )]= \E\left[\exp\left(\imath \la \Xi_{\alpha}, (a^\alpha+b^\alpha )^{\frac{1}{\alpha}} f\ra \right)\right]\eeq{*}
\protect\cite[Definition~2.1]{KumMan}.  Using the above characteristic function~\eqref{eq:limitfield}, it is immediate that the limiting field satisfies this form of stability. An equivalent classical definition, as can be found in \citet[Section~4.8]{LindeBook}, matches ours by means of the Laplacian operator which we introduce as follows.
Choose $a\in \R$. Let us define the operator $(-\Delta)^a$ acting on $L^2(\T^d)$-functions $u$ with Fourier series $\sum_{\nu\in \Z^d}\widehat u(\nu)\mathbf e_{\nu}(\cdot)$ as follows ($(\mathbf e_\nu)_{\nu\in \Z^d}$ denotes a mean-zero orthonormal basis of $L^2(\T^d)$):
\[
(-\Delta)^a \left(\sum_{\nu\in \Z^d}\widehat u(\nu)\mathbf e_{\nu}\right)(\vartheta)=\sum_{\nu\in \Z^d\setminus\{0\}}\|\nu\|^{2a}\widehat u(\nu)\mathbf e_{\nu}(\vartheta).
\]
With this notation we can say the characteristic functional of $\Xi_\alpha$ can be represented as
\begin{equation}\label{eq:limitfield-2}
\E[\exp(\imath \la \Xi_{\alpha}, f\ra )]=\exp\left(-\|(-\Delta)^{-1}f\|_{L^\alpha(\T^d)}^\alpha\right).
\end{equation}
For a reference on $\alpha$-stable cylinder measures one can consult the monograph \citet{LindeBook}. For the reader's convenience, we show that such functionals are well-defined via the Bochner-Minlos theorem (see Appendix~\ref{subsec:app:stable}).    
\begin{remark}
Pluggin in the value $\alpha=2$ in the above Theorem matches the main result of \cite{CHR2016}, concerned specifically with the Gaussian case. 
\end{remark}

The rest of the paper is devoted to the proofs of the above results. 
\section{Proofs on stabilization}\label{sec:stabil}
\subsection{Proof of Lemma~\ref{lemma:infinite mean}}
\noindent
Before we prove the first lemma let us make a general trivial remark. If we assume that $\E[s(o)]=+\infty$ resp. $-\infty$ then necessarily we have that $\E[s^+(o)]= +\infty$ resp. $\E[s^-(o)]=-\infty$ where $s^+$ denotes the positive part and $s^-$ the negative part of the configuration $s$. 
\begin{enumerate}[label=(\roman*),leftmargin=-1pt,rightmargin=-.01pt]
\item By the remark before we can assume that $\E[s^-(o)] < +\infty$, hence $s^-$ is integrable. Note that since the event that $s$ stabilizes is $\Gamma$-invariant, by ergodicity it has probability $0$ or $1$. Assume that $s$ stabilizes almost surely.  
For $M\ge 1$, denote by $$s^M(o):= s(o)\one_{\left\{s(o)\le M\right\}}=s^+(o)\one_{\left \{ 0 \leq s(o)\le M\right\}} - s^-(o)$$ the truncation of the configuration at level $M$.  First we make the following claim:
$$\prob( s^M \text{  stabilizes for all $M\ge 1$})=1.$$
To see this we note that $\mathcal F_s:=\{ f: V\to \R: s+\Delta f\le 1, \,\, f\ge 0\}$ is non-empty if and only if $s$ stabilizes (see \citet[Corollary 2.8]{LMPU}). Now the event that $\mathcal F_{s}\neq \emptyset$ implies that the event $\mathcal F_{s^M}\neq \emptyset$ for all $M\ge 1$, since $s_M\le s$. Hence we have the claim. 

Consequently for any $M$ fixed it holds that $$\E[s^M(o)]=\E[s^+(o)\one_{\left \{ 0 \leq s(o)\le M\right\}}] - \E[s^-(o)]<+\infty$$ and hence applying conservation of density (Proposition 3.1 of \cite{LMPU}) we have that $\E[ s_\infty^M(o)]= \E[ s^M(o)]$. Since the configuration $s_\infty^M$ is stable, $s_\infty^M\le 1$ and so $\E[ s^M(0)]\le 1$ for all $M\ge 1$.  Note that we have on the one side that $s^M(o)$ converges to $s(o)$ almost surely and on the other hand $s^M(o)$ is a monotone increasing sequence in $M$ such that $s^M(o) \geq - s^-(o)$ where $s^-(o)>0$ was assumed integrable. Hence by Fatou's lemma we would get
\[
+\infty = \E \left [\liminf_{M\rightarrow +\infty} s^M(o) \right] \leq \liminf_{M \rightarrow +\infty} \E \left[ s^M(o) \right] \leq 1,
\]
a contradiction.
\item {Since $\E[s(o)]=-\infty$ we can find $M\in (-\infty,\,0]$ such that $\E\left[s(o)\one_{\{s(o)\ge M\}}\right]<1$. Having $s(o)\le s(o)\one_{\{s(o)\ge M\}}$ with probability one, stability follows from \citet[Lemma~4.2]{LMPU}.}
\end{enumerate}
\qed
\subsection{\texorpdfstring{The $\alpha$-singly transient case}{}}\label{subsubsec:singly_trans}
The proof in this case requires a central limit type theorem which we recall here for the reader's convenience.
\begin{theorem}[{Lindeberg-Feller type stable limit theorem, \citet[Theorem~1.1]{DomJun}}.]\label{thm:LF} Suppose
$(\xi_{k,\,j})_{k,\,j\in \N}$ is an i.i.d. array of centered random variables in the domain of
normal attraction of $S\alpha S(1)$, $\alpha\in (0,\,2]$, that is,
\[
\lim_{n\to+\infty}n^{-\nicefrac{1}{\alpha}}\sum_{k=1}^n \xi_{k,\,j}\stackrel{d}{=}S\alpha S (1),\quad \forall\,j\in \N.
\]
Let $\left(u^{(j)}\right)_{j\in \N}$ is a sequence of vectors in $\ell_\alpha$ , i.e. $u^{(j)} := \left(u_k^{(j)}\right)_{k\in \N}\in \ell_\alpha$ for all $j\in \N$. If both
\begin{enumerate}[label=(\arabic*),ref=(\arabic*)]
\item \label{(1)} $\lim_{j\to+\infty}\left\|u^{(j)}\right\|_\alpha=\mathfrak c$,
\item \label{(2)} $\lim_{j\to+\infty}\left\|u^{(j)}\right\|_\infty=0$
\end{enumerate}
hold, then $\sum_k u_k^{(j)}\xi_{k,\,j}<+\infty$ a.~s. for all $j\in \N$ and
\[
\lim_{j\to+\infty}\sum_{k\in \N} u_k^{(j)}\xi_{k,\,j}\stackrel{d}{=}S\alpha S(\mathfrak c).
\]
\end{theorem}
\begin{proof}[Proof of Theorem~\ref{thm:singly_trans}]
We proceed as in \citet[Theorem~3.5]{FMR2009}, \citet[Lemma 5.1]{LMPU}. Assume on the contrary that $s$ stabilizes with probability one. Let $V_1\subset V_2\subset\ldots$ with $\cup_{n\ge 1}V_n=V$. Then using a nested toppling procedure (we stabilize in each volume $V_n$ successively)
\[
s+\Delta u_n=\xi_n, \quad n\in \N
\]
with $\xi_n\le 1$. Let $g_n(x,\,y)$ be the expected number of visits to $y$ by a simple random walk started at $x$ and killed on exiting $V_n$. It holds that \cite[Equation (12)]{LMPU} 
\[
u_n(y)=r^{-1}\sum_{x\in V_n}g_n(x,\,y)(s(x)-1)+r^{-1}\sum_{x\in V_n}g_n(x,\,y)(1-\xi_n(x))
\]
where $r$ is the common degree. Let
\[
\nu_{n,\,\alpha}:=\left(\sum_{y\in V_n}g_n(o,\,y)^\alpha\right)^{\frac{1}{\alpha}}.
\]
We observe that
\eq{}\label{eq:reduction}
P\left(u_n(o)\ge\epsilon \nu_{n,\,\alpha}\right)\ge P\left(\nu_{n,\,\alpha}^{-1}\,\sum_{y\in V_n}g_n(0,\,y)(s(y){{-1}})\ge r\epsilon\right).
\eeq{}
To analyse the right-hand side, we need the following Claim.
\begin{claim}\label{claim:conv_SaS}
\eq{}\label{eq:DOAstable}
\nu_{n,\,\alpha}^{-1}\,\sum_{y\in V_n}g_n(0,\,y)(s(y){{-1}})
\eeq{}
converges in law to a non-degenerate $S\alpha S(1)$ random variable as $n\to+\infty$.
\end{claim}
Since $u_n(o)\uparrow u_\infty(o)$ as $n\to+\infty$ and that we have assumed stabilization, the left-hand side of \eqref{eq:reduction} converges to $0$, while the right-hand side is strictly positive by Claim~\ref{claim:conv_SaS}. This gives a contradiction.
\end{proof}

Let us go into the proof of Claim~\ref{claim:conv_SaS}. 
\begin{proof}[Proof of Claim~\ref{claim:conv_SaS}]
Observe that $s(x)-1$ is a centered random variable for all $x\in V$. Furthermore it belongs to the domain of attraction of an $S\alpha S(1)$.
Let us then verify \ref{(1)}-\ref{(2)} for $\nu_{n,\,\alpha}^{-1}\sum_{y\in V_n}g_n(0,\,y)s(y)$. Taking up the notation of Theorem~\ref{thm:LF}, we define for each $j\in\N$ a sequence $\left(u^{(j)}_k\right)_{k\in\N}$ as follows: if we enumerate the points in $V_j$ such that $V_j=\left\{y_1,\,\ldots,\,y_{|V_j|}\right\}$, let us put
\[
u_k^{(j)}:=
\begin{cases}\nu_{j,\,\alpha}^{-1}\,g_j(o,\,y_k)&k=1,\,\ldots,\,|V_j|\\
0&\text{otherwise}
\end{cases}.
\]
This sequence belongs to $\ell_\alpha$ for fixed $j$ since
\[
\left\|u^{(j)}\right\|_\alpha^\alpha=\sum_{k\in \N}\left(u_k^{(j)}\right)^\alpha=\frac{1}{\nu_{j,\,\alpha}^\alpha}\sum_{y\in V_j}g_j(o,\,y)^\alpha=1.
\]
The above calculation clearly gives that $\lim_{j\to+\infty}\left\|u^{(j)}\right\|_\alpha=1$, so that \ref{(1)} is satisfied. As for \ref{(2)} observe that the boundedness of $g(o,\,\cdot)$ and \eqref{eq:assumption_g1} give
\[
\lim_{j\to+\infty}\frac{g_j(o,\,y)}{\nu_{j,\,\alpha}}=0.
\]
This concludes the proof.
\end{proof}
\subsection{\texorpdfstring{The $\alpha$-doubly transient case}{}}\label{subsubsec:doubly_trans}
In order to characterize the behavior of the divisible sandpile in the $\alpha$ doubly transient case, 
we rely on a result inspired by \citet[Lemma~5.5]{LMPU}, and hence we will postpone its proof to the Appendix in Section~\ref{sec:app:Lemma}.
\begin{lemma}\label{lemma:5_points}
Let $\{y_i\}_{i\ge1}$ be an enumeration of the group $G$. For $\gamma\in \Gamma$, $x\in V$ and $r$ the common degree on $V$ we define
\begin{equation}\label{eq:def:v_e}
v_\gamma(x):=\frac1{r} \sum_{i=1}^\infty g(x, \gamma y_i) Y_{\gamma y_i}.
\end{equation}
Let $e$ be the identity element of $\Gamma$. Then
\begin{enumerate}[label=(\Roman*)]
\item $v_e(o)$ is convergent almost surely.
\item $v_e(o)$ is $\Gamma$-invariant.
\item $v_e(o)$ is almost surely unbounded below.
\end{enumerate}
\end{lemma}
We are now ready to show the main result for the doubly transient case.
\begin{proof}[Proof of Theorem~\ref{thm:doubly_trans}]
Suppose $s$ stabilizes almost surely with odometer $u_\infty$. Then as in Lemma 5.5 of \cite{LMPU} we have that $v_e$ (defined in~\eqref{eq:def:v_e}) has $\Gamma$-invariant law and $\Delta v_e=1-s$, hence $h= v-u_\infty$ has invariant law and is harmonic on $V$. Observe that the assumptions of \citet[Theorem~2.2]{KokTaq} are satisfied, in such a way that we can conclude that
$$\lim_{t\to+\infty}\frac{\prob(\left| v_e(o)\right|> t)}{\prob(\left|Y_o\right|>t)}=\frac{1}{r}\sum_{i=1}^{+\infty} g(o,\, y_i)^\alpha .$$
As a consequence $v_e$ has a right regularly varying tail of index $-\alpha$ and hence $$\E\left[ \left|v_e(o)\right|^{\alpha-\eps}\right]<+\infty$$
for all $0<\eps<\alpha$. Hence by \citet[Lemma~5.4]{LMPU}, we have that $h$ is constant almost surely. Since $u_\infty\ge 0$ and $v_e$ is unbounded below almost surely by Lemma~\ref{lemma:5_points}, we have a contradiction.
\end{proof}

\section{Proof of Theorem~\ref{thm:scaling_RV}}\label{sec:odometer_limit}
\subsection{Preliminaries}
Consider the Hilbert space $L^2(\Z_n^d)$ of complex valued functions on the discrete torus endowed with the inner product
$$\la f, g \ra = \frac1{n^d} \sum_{x\in \Z_n^d} f(x)\overline{g(x)}.$$
The Pontryagin dual group of $\Z_n^d$ is identified again with $\Z_n^d$. Let $\{ \chi_w: w\in \Z_n^d\}$ denote the characters of the group where $\chi_w(x)= \exp(2\pi \mathbbm{i} x\cdot {w}/{n})$. 
The eigenvalues of the Laplacian $\Delta_g$ on discrete tori are given by $$\lambda_w= -4 \sum_{i=1}^d \sin^2\left(\frac{\pi w_i}{n}\right),\quad w\in\Z_n^d.$$
We use the shortcut $g_x(y):= g(y,x)$.
Let $\widehat g_x$ denote the Fourier transform of $g_x$.  
\subsection{Proof of Theorem~\ref{thm:scaling_RV}}
The proof of Theorem~~\ref{thm:scaling_RV} relies on two steps. As done in \cite{CHR2016}, the proof is based on determining the scaling limit in a ``simpler'' case, that is, when the variables $\sigma$ in Proposition~\ref{prop:1.3LMPU} are i.i.d.~$S\alpha S(1)$. Then we will conclude in the more general case using the theorem for symmetric stable laws.
\subsubsection{Proof for the $\alpha$-stable case}
By means of Proposition~\ref{prop:1.3LMPU} and the fact that all test functions have mean $0$, the main Theorem on the scaling limit of the odometer will follow once we prove this statement:
\begin{theorem}\label{thm:scaling_RV_stable}
Let $(\sigma(x))_{x\in \Z_n^d}$ be i.i.d. $S\alpha S(1)$ random variables. For all $f\in C^\infty(\T^d)$ with mean zero, the variables $\la \Xi_n, f \ra$ converges in distribution to $\la \Xi_{\alpha}, f\ra$ where $\Xi_\alpha$ is the same of Theorem~\ref{thm:scaling_RV}.
\end{theorem}
\begin{proof}[Overview of the proof]
Let us denote by $v_n(y)=(2d)^{-1} \sum_{x\in \Z_n^d} g(x,y)(s(x)-1)$ and as $u(\cdot)$ the odometer function. Note that it follows from  Proposition~\ref{prop:1.3LMPU} that the odometer has the following representation: 
\begin{equation}
u(x) = v_n(x)-\min_{z\in \Z_n^d}v_n(z).
\end{equation}
Let us define the following functional: for any function $h_n: \Z_n^d\to \R$ set
$$\Xi_{h_n}(x):={4\pi^2}\sum_{z\in \T^d_n}n^{d-\frac{d}{\alpha}-2}h_n({nz})\one_{B\left(z,\,\nicefrac{1}{2n}\right)}(x),\quad x\in \T^d.$$
For $f\in C^\infty(\T^d)$ such that $\int_{\T^d} f(x)\De x=0$ it follows immediately that
$$\la \Xi_{u}, f\ra = \la \Xi_{v_n}, f\ra.$$
If we call 
\[w_n(y):=\left(2d\right)^{-1}\sum_{x\in \Z_n^d} g(x,y)\sigma(x),\]
by the mean-zero property of the test functions and the Random Target Lemma (see Section 5 of \cite{CHR2016}) we deduce that $\la \Xi_{v_n}, f\ra = \la \Xi_{w_n}, f \ra.$ Therefore we shall reduce ourselves to study the convergence of the field $\Xi_{w_n}$.

The proof consists of 5 steps, which we will elucidate here together with some notation. Later we will show each step separately.
We write $c_n:={4\pi^2} n^{d-\nicefrac{d}{\alpha}-2}$. Let us denote by 
\eq{}\label{eq:def_H}H_n(z)= \int_{B\left(z,\,\frac{1}{2n}\right)} f(t)\De t.\eeq{}
We can then rewrite
\begin{align}
\la \Xi_{w_n}, f \ra &= c_n\sum_{z\in \T_n^d} w(nz)H_n(z) \nonumber\\
&= \sum_{x\in \Z_n^d} \left( c_n (2d)^{-1}\sum_{z\in \T_n^d} g(x,nz) H_n(z) \right) \sigma(x)=\sum_{x\in \Z_n^d} k_n(x)\sigma(x),\label{eq:for_later}
\end{align}
where 
\eq{}\label{eq:def_k_n}
k_n(x):= c_n(2d)^{-1}\sum_{z\in T_n^d} g(x,nz) H_n(z),\quad x\in \Z_n^d .\eeq{}
Hence using the characteristic function of $\alpha$-stable variables
\[
\E[ \exp( \imath \la \Xi_{w_n}, f \ra)]=\exp\left(-\sum_{x\in \Z_n^d} |k_n(x)|^\alpha\right).
\]
Letting $L_n(z):= H_n(z/n)$, we rewrite (using Perseval's lemma)
\begin{align}
k_n(x)&=c_n (2d)^{-1}\sum_{z\in \T_n^d} g(x,nz) H_n(z) = c_n (2d)^{-1} \sum_{z\in \Z_n^d} g(x,z)L_n(z) \,\, \,  \nonumber\\
&= c_n(2d)^{-1}n^d\la g_x, L_n \ra=  c_n  (2d)^{-1}n^d\sum_{z\in \Z_n^d} \widehat{g_x}(z) \widehat{L_n}(z)\label{def:L_n}
\end{align}
for $x\in\Z_n^d$. Now we will split the above sum into contributions from the site $z=0$ and from other sites. Note that $\widehat{g_x}(0)$ is independent of $x$ (cf. Equation (3.1) of \citet{CHR2016}). 
Moreover
\begin{align*}
\widehat L_n(0)&=  n^{-d} \sum_{z\in \Z_n^d} L_n(z) =n^{-d}\sum_{z\in \T_n^d} H_n(z)\\
&=n^{-d}\sum_{z\in \T_n^d} \int_{B\left(z, \frac{1}{2n}\right)} f(u)\De u =n^{-d} \int_{\T^d} f(u)\De u=0. 
\end{align*}
We can use the fact that \cite[Equation (20)]{LMPU}
\eq{}\label{eq:fundamental}
\lambda_a \widehat{g_x}(a)=-2d n^{-d}\chi_{-a}(x),\quad a\neq 0
\eeq{}
to deduce that
\begin{align*}
c_n n^d& (2d)^{-1}\sum_{z\in \Z_n^d\setminus\{0\}} \widehat{g_x}(z) \widehat{L_n}(z)=-c_n \sum_{z\in \Z_n^d\setminus\{0\}} \frac{ \chi_{-z}(x)}{\lambda_z}\widehat{L_n}(z)\\
&= -c_n\sum_{z\in \Z_n^d\setminus\{0\}} \frac{ \chi_{-z}(x)}{\lambda_z} \la L_n, \chi_z\ra= -\frac{c_n}{n^d} \sum_{z\in \Z_n^d\setminus\{0\}} \frac{ \chi_{-z}(x)}{\lambda_z} \sum_{w\in \Z_n^d} L_n(w) \chi_{-z}(w)\\
&=-\frac{c_n}{n^d}\sum_{w\in\T_n^d}\int_{B\left(w,\,\frac{1}{2n}\right)} f(u)\De u \sum_{z\in \Z_n^d\setminus\{0\}} \frac{ \chi_{-z}(x)\chi_{-z}(nw)}{\lambda_z} .
\end{align*}
{Defining 
\[R_n(w):= \int_{B\left(w,\,\frac{1}{2n}\right)} (f(u)-f(w))\De u\]} we can split further the integral in the above equality and obtain
\begin{align}
k_n(x)&=-\frac{c_n}{n^d}\sum_{w\in\T_n^d}\left( n^{-d} f(w)+R_n(w)\right) \sum_{z\in \Z_n^d\setminus\{0\}} \frac{ \chi_{-z}(x)\chi_{-z}(nw)}{\lambda_z}\nonumber\\
&=-\frac{c_n}{n^{2d}}\sum_{w\in\T_n^d} f(w)\!\sum_{z\in \Z_n^d\setminus\{0\}} \frac{ \chi_{-z}(x)\chi_{-z}(nw)}{\lambda_z}-\frac{c_n}{n^d}\sum_{w\in\T_n^d}R_n(w)\!\sum_{z\in \Z_n^d\setminus\{0\}} \frac{ \chi_{-z}(x)\chi_{-z}(nw)}{\lambda_z}\nonumber\\
&:= l_n(x)+ C_n(x)\label{eq:def_lC}.
\end{align}
Now our first step is to show that the convergence of $\exp\left(-\sum_{x\in \Z_n^d} |k_n(x)|^\alpha\right)$ can be given in terms of the same quantity where $k_n(\cdot)$ is replaced by $l_n(\cdot)$:
\begin{step}\label{step:1}
\[
\lim_{n\to+\infty}\left|\exp\left(-\sum_{x\in \Z_n^d} |k_n(x)|^\alpha\right)- \exp\left(-\sum_{x\in \Z_n^d} |l_n(x)|^\alpha\right)\right|=0.
\]
\end{step}
The next steps aim at proving that $l_n$ is giving us the correct characteristic function. {In Step~\ref{step:2} we are introducing a mollifier which will help to extend sums from $\Z_n^d$ to the whole lattice.} 
\begin{step}\label{step:2}
{Let $\phi\in \mathcal S(\R^d)$, the Schwartz space, with
$\int_{\R^d}\phi(x)\De x=1$. Let $\eps>0$ and let $\phi_\eps(x):= \eps^{-d}\phi\left({x}{\eps}^{-1}\right)$ for $\eps>0$.} Then 
\begin{align*}
\lim_{\eps\downarrow 0}\lim_{{n\to+\infty}}\left|\sum_{x\in \Z_n^d} |l_n(x)|^{\alpha}- \frac{c_n^\alpha}{n^{d\alpha}} \sum_{x\in \T_n^d}\left| \sum_{z\in \Z_n^d\setminus \{0\}} \frac{\widehat \phi_\eps(z)\exp(-2\pi \imath z\cdot x) }{\lambda_z} \widehat{f_n}(z)\right|^{\alpha}\right|=\mathrm{O}\left({\eps^{\min\{\alpha,\,1\}}}\right)
\end{align*}
where $\widehat{f_n}(z)= n^{-d}\sum_{w\in \T_n^d} f(w)\exp(2\pi \mathbbm{i} w\cdot z)$. 
\end{step}
The goal of the third step is to approximate each eigenvalue $\lambda_z$ of the Laplacian with the norm of the point $z$, namely
\begin{step}\label{step:3}
For all $\eps>0$
\begin{align*}
\lim_{n\to+\infty}\frac{c_n^\alpha}{n^{d\alpha}} &\left|\sum_{x\in \T_n^d}\left[\left| \sum_{z\in \Z_n^d\setminus \{0\}} \frac{\widehat \phi_\eps(z)\e^{-2\pi \imath z\cdot x} }{\lambda_z} \widehat{f_n}(z)\right|^{\alpha}
\right.\right.\\
&\left.\left.-\frac{ n^{2\alpha}}{4^\alpha {\pi^{2\alpha}} } \left| \sum_{z\in \Z_n^d\setminus \{0\}} \frac{\widehat \phi_\eps(z)\e^{-2\pi \imath z\cdot x} }{\|z\|^2} \widehat{f_n}(z)\right|^{\alpha}\right]\right|=0
\end{align*}
\end{step}
In the next step we extend the sums in Step \ref{step:3} over $\Z^d$ using the decay of the mollifier.
\begin{step}\label{step:4}
For all $\eps>0$
\begin{align*}
\lim_{n\to+\infty}\frac{c_n^\alpha n^{2\alpha}}{n^{d\alpha}{4^\alpha\pi^{2\alpha}}}&\left| \sum_{x\in \T_n^d}\left| \sum_{z\in \Z_n^d\setminus \{0\}} \frac{\widehat \phi_\eps(z)\exp(-2\pi \imath z\cdot x) }{\|z\|^2} \widehat{f_n}(z)\right|^{\alpha}\right.\\
& -
\left. \sum_{x\in \T_n^d}\left| \sum_{z\in \Z^d\setminus \{0\}} \frac{\widehat \phi_\eps(z)\exp(-2\pi \imath z\cdot x) }{\|z\|^2} \widehat{f_n}(z)\right|^{\alpha}\right|=0.
\end{align*}
\end{step}
At last, we can finally show the convergence of the sum to the required integral.
\begin{step}\label{step:5}
\begin{align*}
\lim_{\eps\downarrow 0}\lim_{{n\to+\infty}}\frac{c_n^\alpha n^{2\alpha}}{n^{d\alpha}{4^\alpha\pi^{2\alpha}}}\sum_{x\in \T_n^d}&\left| \sum_{z\in \Z^d\setminus \{0\}} \frac{\widehat \phi_\eps(z)\exp(-2\pi \imath z\cdot x) }{\|z\|^2} {\widehat{f_n}}(z)\right|^{\alpha}=\\
& \int_{\T^d} \left| \sum_{z\in \Z^d\setminus \{0\}} \frac{\exp(-2\pi \imath z\cdot x)}{\|z\|^2} \widehat{f}(z)\right|^{\alpha} \De x.
\end{align*}
\end{step}
\end{proof}
{The core of the proof is showing the 5 steps. They are logically dependent one from another as follows:}
\begin{center}
{Step~\ref{step:5}$\Rightarrow$ Step~\ref{step:4}$\Rightarrow$ Step~\ref{step:3}$\Rightarrow$
Step~\ref{step:2}$\Rightarrow$
Step~\ref{step:1}.}
\end{center}
{We will now begin to show the proof of each step assuming the subsequent ones, and will finally conclude with Step~\ref{step:5}.}
\begin{proof}[Proof of Step \ref{step:1}]
Let us denote by
\eq{}\label{eq:def_t_n}
t_n(x):=t k_n(x)+(1-t)l_n(x),\quad t\in [0,\,1].
\eeq{}
 Using $|\exp(-a)-\exp(-b)|\le |a-b|$ for $a,\, b\ge 0$ we get
\begin{align}
&\left| \exp\left(-\sum_{x\in \Z_n^d} |k_n(x)|^\alpha\right)- \exp\left(-\sum_{x\in \Z_n^d} |l_n(x)|^\alpha\right)\right|\nonumber\\
&\le \left| \sum_{x\in \Z_n^d} |k_n(x)|^\alpha-|l_n(x)|^\alpha\right|\le \sum_{x\in \Z_n^d} \left| |k_n(x)|^\alpha-|l_n(x)|^\alpha\right|\label{eq:bath} .
\end{align}
By the mean value theorem we can bound the last term as follows:
\eq{}
\begin{cases}
\sum_{x\in \Z_n^d} \alpha |t_n(x)|^{\alpha-1} |C_n(x)|\label{eq:Holder_first} &\text{if }\alpha>1 \\
\sum_{x\in \Z_n^d} | C_n(x)|^\alpha &\text{if }\alpha\le 1
\end{cases}.
\eeq{}
From \eqref{eq:def_t_n}, \eqref{eq:Holder_first} and the bound 
\[
(a+b)^r\le 2^{r}\left(a^r+b^r\right),\quad a,\,b\ge 0,\,r\ge 0
\]
we get
\begin{align}\label{eq:fundamental_2}
\sum_{x\in \Z_n^d}& \left| |k_n(x)|^\alpha-|l_n(x)|^\alpha\right|\nonumber\\
&\le\begin{cases}
{\alpha 2^{\alpha-1} \sum_{x\in \Z_n^d} \left|C_n(x) \right|^{\alpha} + \alpha 2^{\alpha-1} \sum_{x\in \Z_n^d} \left|C_n(x) \right| \left|l_n(x)\right|^{\alpha-1}}&\text{if }\alpha>1\\
\sum_{x\in \Z_n^d} | C_n(x)|^\alpha &\text{if }\alpha\le 1
\end{cases}.
\end{align}
Let us look at $\sum_{x\in \Z_n^d}\left|C_n(x)\right|^\alpha$. We notice that 
\eq{}\label{eq:before_holder}\left(\sum_{x\in \Z_n^d}\left|C_n(x)\right|^\alpha\right)^{1/\alpha}= n^{\frac{d}{\alpha}}\left(n^{-d}\sum_{x\in \Z_n^d}\left|C_n(x)\right|^\alpha\right)^{1/\alpha}
\eeq{}
Observe that by H\"older's inequality we have that 
$$
\left(n^{-d}\sum_{x\in \Z_n^d}\left|C_n(x)\right|^\alpha\right)^{1/\alpha} \le\left(n^{-d}\sum_{x\in \Z_n^d}\left|C_n(x)\right|^2\right)^{\frac12}=\|C_n\|_2.$$ 
Hence an appropriate bound on the $L_2$-norm of $C_n$ will suffice to prove that this term is small.
First we provide a crude bound for $C_n(x)$:
\begin{align}
C_n(x) &=\frac{c_n}{n^d}\sum_{w\in\T_n^d}R_n(w)\sum_{z\in \Z_n^d\setminus\{0\}} \frac{ \chi_{-z}(x)\chi_{-z}(nw)}{\lambda_z}\nonumber\\
&=c_n \sum_{z\in \Z_n^d\setminus\{0\}}\frac{ \chi_{-z}(x)}{\lambda_z} n^{-d} \sum_{w\in \Z_n^d} R_n(w/n) \chi_{-z}(w)=4 \pi^2 n^{d-\frac{d}{\alpha}-2}\sum_{z\in \Z_n^d\setminus\{0\}}\frac{ \chi_{-z}(x)}{\lambda_z}\widehat{\mathcal R_n}(z)\nonumber\\
&= 4\pi^2 n^{-\frac{d}{\alpha}-2}\sum_{z\in \Z_n^d\setminus\{0\}}\frac{ \chi_{-z}(x)}{\lambda_z}n^d\widehat{\mathcal R_n}(z)\label{eq:crude_C}
\end{align}
where $\mathcal R_n(w):= R_n(w/n)$. We wish to bound the $L_2$-norm of $C_n$ and to do so we employ \citet[Lemma 7]{CHR2016}.
It follows from it and \eqref{eq:crude_C} that
\begin{align*}
\|C_n\|_2^2 &=(4\pi^2)^2 n^{2\left(d-\frac{d}{\alpha}-2\right)}\sum_{z\in \Z_n^d\setminus \{0\}}\frac{ \left|\widehat{\mathcal R_n}(z)\right|^2}{|\lambda_z|^2}\le C n^{2\left(d-\frac{d}{\alpha}\right)}\sum_{z\in \Z_n^d}\left|\widehat{\mathcal R_n}(z)\right|^2= C n^{d-\frac{2d}{\alpha}}\sum_{z\in \Z_n^d}\left|\mathcal R_n(z)\right|^2\\
&= C n^{d-\frac{2d}{\alpha}}\sum_{z\in \T_n^d}\left| R_n(z)\right|^2\le C n^{-\frac{2d}{\alpha}-{2}}.
\end{align*}
Note that in the last step we have used that
$$|R_n(w)|\le \int_{B\left(w,\,\frac1{2n}\right)} |f(u)-f(w)| \De u\le \|\nabla f\|_{\infty} n^{-d-1}.$$
We have deduced that
\eq{}\label{eq:norm_2} \|C_n\|_2\le  C n^{-\frac{d}{\alpha}-{1}}.\eeq{}
This plugged into \eqref{eq:before_holder} shows that the first summand of the first line resp. the second line of \eqref{eq:fundamental_2} tends to zero.

As for the second summand of the first line in \eqref{eq:fundamental_2}, we wish to apply H\"older's inequality:
\begin{align}
\sum_{x\in \Z_n^d} |C_n(x)| &|l_n(x)|^{\alpha-1}\le n^d \left( n^{-d} \sum_{x\in \Z_n^d} |C_n(x)|^{\alpha}\right)^{\frac{1}{\alpha}}\left( n^{-d} \sum_{x\in \Z_n^d} \left(|l_n(x)|^{\alpha-1}\right)^{\frac{\alpha}{\alpha-1}}\right)^{\frac{\alpha-1}{\alpha}}\nonumber\\
&\overset{\alpha\le 2}\le n^d\left( n^{-d} \sum_{x\in \Z_n^d} |C_n(x)|^{2}\right)^{\frac{1}{2}}\left( n^{-d} \sum_{x\in \Z_n^d} |l_n(x)|^{\alpha}\right)^{\frac{\alpha-1}{\alpha}}\nonumber\\
&= n^{d} \|C_n\|_2 n^{-\frac{d(\alpha-1)}{\alpha}} \left( \sum_{x\in \Z_n^d} |l_n(x)|^{\alpha}\right)^{\frac{\alpha-1}{\alpha}}\nonumber\\
&\overset{\eqref{eq:norm_2}} \le n^d n^{-\frac{d}{\alpha}-1}n^{-\frac{d(\alpha-1)}{\alpha}} \left( \sum_{x\in \Z_n^d} |l_n(x)|^{\alpha}\right)^{\frac{\alpha-1}{\alpha}}=n^{-1}\left( \sum_{x\in \Z_n^d} |l_n(x)|^{\alpha}\right)^{\frac{\alpha-1}{\alpha}}\label{eq:number13}
\end{align}
Now in Steps \ref{step:3}-\ref{step:4}-\ref{step:5} we shall show that 
$$\lim_{n\to+\infty}\sum_{x\in \Z_n^d} |l_n(x)|^{\alpha}= \int_{\T^d} \left| \sum_{z\in \Z^d\setminus \{0\}} \frac{\exp(-2\pi \imath z\cdot x)}{\|z\|^2} \widehat{f}(z)\right|^{\alpha} \De x.$$
Hence \eqref{eq:number13} and consequently the second summand in the first inequality of \eqref{eq:fundamental_2} tends to zero. This concludes the proof of the first step.
\end{proof}
\begin{proof}[Proof of Step \ref{step:2}]
Recall that we have
\begin{align*}
\sum_{x\in \Z_n^d} |l_n(x)|^{\alpha}&=  \sum_{x\in \T_n^d} \left|\frac{c_n}{n^{d}} \sum_{w\in \T_n^d} n^{-d} f(w)\sum_{z\in \Z_n^d\setminus \{0\}} \frac{\exp(-2\pi \imath z\cdot x/n) \exp(2\pi \imath z\cdot w)}{\lambda_z}\right|^\alpha.
\end{align*}
Let us write as before $l_n(x)$ as sum of two quantities:
\begin{align*}
l_n(x)
&=\frac{c_n}{n^d}\sum_{w\in \T_n^d} n^{-d} f(w)\sum_{z\in \Z_n^d\setminus \{0\}} \left(1-\widehat{\phi_\eps}(z)\right)\frac{\exp(-2\pi \imath z\cdot x/n) \exp(2\pi \imath z\cdot w)}{\lambda_z}+\\
&+\frac{c_n}{n^d}\sum_{w\in \T_n^d} n^{-d} f(w)\sum_{z\in \Z_n^d\setminus \{0\}} \widehat{\phi_\eps}(z)\frac{\exp(-2\pi \imath z\cdot x/n) \exp(2\pi \imath z\cdot w)}{\lambda_z}\\
&=: C_n^{(1)}(x)+ l_n^{(1)}(x).
\end{align*}
Exactly as in \eqref{eq:Holder_first} one has
\begin{align}\label{eq:nepal}
&\left|\sum_{x\in \Z_n^d} |l_n(x)|^{\alpha}-\sum_{x\in \Z_n^d} \left|l_n^{(1)}(x)\right|^{\alpha}\right|\nonumber\\
&\le
\begin{cases}
 \alpha 2^{\alpha-1}\sum_{x\in \Z_n^d}\left| C_n^{(1)}(x)\right|^\alpha+ \alpha 2^{\alpha-1}\sum_{x\in \Z_n^d} \left|C_n^{(1)}(x)\right| \left|l_n^{(1)}(x)\right|^{\alpha-1}&\text{if }\alpha> 1\\
 \sum_{x\in\Z_n^d}\left|C_n^{(1)}(x)\right|^\alpha&\text{if }\alpha\le 1
 \end{cases}.
\end{align}
As before in Step \ref{step:1}, we show the terms on the right-hand side go to zero. Let us look at the first sum. 
$$\left(\sum_{x\in \Z_n^d}\left|C_n^{(1)}(x)\right|^\alpha\right)^{\frac1\alpha}= n^{\frac{d}{\alpha}}\left(n^{-d}\sum_{x\in \Z_n^d}\left|C_n^{(1)}(x)\right|^\alpha\right)^{\frac1\alpha}= n^{\frac{d}{\alpha}}\| C_n^{(1)}\|_\alpha.$$
Observe that by H\"older's inequality, using $\alpha<2$, we have that $n^{\frac{d}{\alpha}}\left\| C_n^{(1)}\right\|_\alpha\le n^{\frac{d}{\alpha}}\left\|C_n^{(1)}\right\|_2$. Hence again it all boils down to finding an estimate for $\left\|C_n^{(1)}\right\|_2$.
Recall that
\begin{align*}
C_n^{(1)}(x)= \frac{c_n}{n^d}\sum_{z\in \Z_n^d\setminus \{0\}} \left(1-\widehat{\phi_\eps}(z)\right)\frac{\e^{-2\pi \imath z\cdot \frac{x}{n}}}{\lambda_z}\widehat{f_n}(z).
\end{align*}
Now note that, since $\left|1-\widehat{\phi_\eps}(z)\right|\le C\eps\|z\|$ as proved by \citet[Eq. (2.11)]{CHR2016},
\begin{align*}
&n^{-d} \sum_{x\in \Z_n^d} \left|C_n^{(1)}(x)\right|^2\\
&=n^{-d}\frac{c_n^2}{n^{2d}}\sum_{x\in \Z_n^d}\sum_{z,z'\in \Z_n^d\setminus \{0\}} \left(1-\widehat{\phi_\eps}(z)\right)\left(1-\widehat{\phi_\eps}(z')\right)\frac{\e^{-2\pi \imath z\cdot \frac{x}{n}}\e^{2\pi \imath z'\cdot \frac{x}{n} } }{\lambda_z \lambda_{z'}}\widehat{f_n}(z)\overline{\widehat{f_n}(z')}\\
&=n^{-2d}c_n^2\sum_{z\in \Z_n^d\setminus \{0\}}\frac{ \left|1-\widehat{\phi_\eps}(z)\right|^2|\widehat{f_n}(z)|^2}{|\lambda_z|^2}\le C n^{-2d+4}c_n^2\sum_{z\in \Z_n^d\setminus \{0\}}\frac{ \left|1-\widehat{\phi_\eps}(z)\right|^2|\widehat{f_n}(z)|^2}{\|z\|^4} 
\end{align*}
using the bound of \citet[Lemma 7]{CHR2016}. We can further bound the last member of the inequality from above with
\begin{align*}
 Cn^{-2d+4}c_n^2\sum_{z\in \Z_n^d\setminus \{0\}}\frac{ \eps^2\|z\|^2|\widehat{f_n}(z)|^2}{\|z\|^4}&\le Cn^{-2d+4}\eps^2 c_n^2\sum_{z\in \Z_n^d\setminus \{0\}} \left|\widehat{f_n}(z)\right|^2\\
&\le C\eps^2 n^{-d -2\frac{d}{\alpha}} \sum_{z\in \Z_n^d} |f_n(z)|^2.
\end{align*}
To sum up, for the first summand of \eqref{eq:nepal} we have obtained a bound of the form
\eq{}\label{eq:norm_C_n}\|C_n^{(1)}\|_2\le n^{-\frac{d}{\alpha}} \eps \left(\frac{1}{n^d}\sum_{z\in \T_n^d} |f(z)|^2\right)^{\frac12}.\eeq{}
Hence for the first term we have
$$n^{\frac{d}{\alpha}}\| C_n^{(1)}\|_\alpha\le \eps\left(\frac{1}{n^d}\sum_{z\in \T_n^d} |f(z)|^2\right)^{\frac12}.$$
Observing that $n^{-d}\sum_{z\in \T_n^d} |f(z)|^2\to\int_{\T^d} |f(z)|^2\De z$ we get the result.

It is time now to handle the second term appearing in \eqref{eq:nepal}. Using H\"older we have that
\begin{align*}
&\sum_{x\in \Z_n^d} |C_n^{(1)}(x)| |l_n^{(1)}(x)|^{\alpha-1}\le n^d \|C_n^{(1)}\|_2 n^{-\frac{d(\alpha-1)}{\alpha}}\left( \sum_{x\in \Z_n^d} |l_n^{(1)}(x)|^{\alpha}\right)^{\frac{\alpha-1}{\alpha}}\\
&\stackrel{\eqref{eq:norm_C_n}}{\le}  \eps \left( \sum_{x\in \Z_n^d} |l_n^{(1)}(x)|^{\alpha}\right)^{\frac{\alpha-1}{\alpha}}\left(\frac{1}{n^d}\sum_{z\in \T_n^d} |f(z)|^2\right)^{\frac12}\\
&=\eps \left( \sum_{x\in \Z_n^d} |l_n^{(1)}(x)|^{\alpha}\right)^{\frac{\alpha-1}{\alpha}}\left(\frac{1}{n^d}\sum_{z\in \T_n^d} |f(z)|^2\right)^{\frac12}.
\end{align*}
Steps \ref{step:3}-\ref{step:4}-\ref{step:5} will show that $\sum_{x\in \Z_n^d} |l_n^{(1)}(x)|^{\alpha}$ converges as $n\to+\infty$ to a finite quantity, and hence the above product will be neglibile in the limit.

\end{proof}
\begin{proof}[Proof of Step \ref{step:3}]
We rewrite
\begin{align*}
&\frac{c_n^\alpha}{n^{d\alpha}} \sum_{x\in \Z_n^d}\left| \sum_{z\in \Z_n^d\setminus \{0\}} \frac{\widehat \phi_\eps(z)\e^{-2\pi \imath z\cdot \frac{x}{n}} }{\lambda_z} \widehat{f_n}(z)\right|^{\alpha}={c_n^\alpha n^{\alpha(2-d)}} \sum_{x\in \Z_n^d}\left| \sum_{z\in \Z_n^d\setminus \{0\}} \frac{\widehat \phi_\eps(z)\e^{-2\pi \imath z\cdot \frac{x}{n}} }{4\|\pi z\|^2} \widehat{f_n}(z)\right.\\
&\left.+\sum_{z\in \Z_n^d\setminus \{0\}} \widehat \phi_\eps(z)\e^{-2\pi \imath z\cdot \frac{x}{n}} \left(\frac1{n^2\lambda_z}-\frac{1}{4\|\pi z\|^2}\right) \widehat{f_n}(z)\right|^{\alpha}=: \sum_{x\in \Z_n^d}\left| l_n^{(2)}(x)+ C_n^{(2)}(x)\right|^{\alpha}.
\end{align*}
{We will only deal here with the case $\alpha>1$. The same procedure of Steps \ref{step:1}-\ref{step:2} can be followed to treat the case $\alpha\le 1$}. We observe that
\begin{align}\frac{c_n^\alpha}{n^{d\alpha}} \sum_{x\in \Z_n^d}\left| \sum_{z\in \Z_n^d\setminus \{0\}} \frac{\widehat \phi_\eps(z)\e^{-2\pi \imath z\cdot \frac{x}{n}} }{\lambda_z} \widehat{f_n}(z)\right|^{\alpha}
&\le \alpha 2^{\alpha-1}\sum_{x\in \Z_n^d}\left|C_n^{(2)}(x)\right|^{\alpha} \nonumber\\
&+\alpha 2^{\alpha-1} \sum_{x\in \Z_n^d} \left|C_n^{(2)}(x)\right| \left| l_n^{(2)}(x)\right|^{\alpha-1}.\label{eq:bath_1}\end{align}
In order to show that the first term goes to zero, it is enough to show that $n^{\frac{d}{\alpha}} \|C_n^{(2)}\|_2$ tends to $0$.
We get
$$C_n^{(2)}(x)= {c_n n^{2-d}} \sum_{z\in \Z_n^d\setminus\{0\}} \widehat \phi_\eps(z) \e^{-2\pi \imath z\cdot \frac{ x}{n}} \left(\frac{1}{n^2\lambda_z}-\frac1{4\|\pi z\|^2}\right)\widehat{f_n}(z).$$
In the same fashion as before, we use the orthogonality of the characters, \citet[Lemma~7]{CHR2016}, the uniform bound on $\left\|\widehat{\phi}\right\|_\infty$ and Parseval's identity to get
\begin{align*}
n^{-d}&\sum_{x\in \Z_n^d} \left|C_n^{(2)}(x)\right|^2 = {c_n^{2} n^{4-2d}} \sum_{z\in \Z_n^d\setminus\{0\}} \left|\widehat \phi_\eps(z)\right|^2\left|\frac{1}{n^2\lambda_z}-\frac1{4\|\pi z\|^2}\right|^2\left|\widehat{f_n}(z)\right|^2\\
&\le C{c_n^2 n^{4-2d}} n^{-4} \sum_{z\in \Z_n^d\setminus\{0\}}\left|\widehat \phi_\eps(z)\right|^2\left|\widehat f_n(z)\right|^2\\
&= C n^{2d-2\frac{d}{\alpha}-4}n^{-2d} \sum_{z\in \Z_n^d}\left|\widehat{f_n}(z)\right|^2\le n^{-\frac{2d}{\alpha}-4} \left(n^{-d}\sum_{z\in \T_n^d}|f(z)|^2\right).
\end{align*}
Hence we have that 
$$\left\|C_n^{(2)}\right\|_2 \le n^{-\frac{d}{\alpha}-2}  \left(n^{-d}\sum_{z\in \T_n^d}|f(z)|^2\right)^{\frac12}$$
showing that $n^{\frac{d}{\alpha}} \|C_n^{(2)}\|_2\to 0$. Now provided we can show Step \ref{step:4} and Step \ref{step:5}, the second term of \eqref{eq:bath_1} would converge to zero  along the lines of \eqref{eq:number13}, completing thus the proof of Step~\ref{step:3}.
\end{proof}
\begin{proof}[Proof of Step \ref{step:4}]
As before we write $l_n^{(2)}(x):= l_n^3(x)+C_n^{(3)}(x)$ where we recall
$$l_n^{(2)}(x)= {c_n n^{2-d}} \sum_{z\in \Z^d\setminus \{0\}} \frac{\widehat \phi_\eps(z)\exp(-2\pi \imath z\cdot x) }{4\|\pi z\|^2} \widehat{f_n}(z)$$
and set
$$C_n^{(3)}(x):= {c_n n^{2-d}} \sum_{\|z\|_\infty>n} \frac{\widehat \phi_\eps(z)\exp(-2\pi \imath z\cdot x) }{4\|\pi z\|^2} \widehat{f_n}(z).$$
We now show that $n^{\frac{d}{\alpha}} \|C_n^{(3)}\|_2$ tends to $0$. Using orthogonality and the approximation of Euler-MacLaurin's formula \cite[Theorem 1]{Apostol} we get that
\begin{align*}
n^{-d}&\sum_{x\in \Z_n^d} \left|C_n^{(3)}(x)\right|^2=\frac{c_n^2 n^{4-2d}}{16}\sum_{\|z\|_\infty>n} \frac{|\widehat \phi_\eps(z)|^2}{\|\pi z\|^4} \left|\widehat{f_n}(z)\right|^2\le \|f\|_\infty^2{n^{-2\frac{d}{\alpha}}}\sum_{\|z\|_\infty>n} \frac{|\widehat \phi_\eps(z)|^2}{\|z\|^4} \\
&\le  \|f\|_\infty^2{n^{-2\frac{d}{\alpha}}}\sum_{\|z\|_\infty>n} \frac{1}{\|z\|^4(1+\|z\|)^{d+1}}\\
&\le\|f\|_\infty^2{n^{-2\frac{d}{\alpha}}}\int_n^{+\infty} t^{d-1} t^{-d-5}\De t+C{n^{-2\frac{d}{\alpha}-6}}\le C n^{-\frac{2d}{\alpha}-5}.
\end{align*}
We have used here that ${|\widehat{f_n}(z)|\le \|f\|_\infty}$ and the fast decay of $\widehat{\phi_\eps}$ at infinity. Hence we have that $n^{\frac{d}{\alpha}} \|C_n^{(3)}\|_2\le Cn^{-5/2}$. Since the conclusion follows similarly to Steps~\ref{step:1}-\ref{step:2} we skip the rest of the proof. 
\end{proof}
\begin{proof}[Proof of Step \ref{step:5}]
By our choice of $c_n$ we have
$$\frac{c_n^\alpha n^{2\alpha+d}}{(4\pi^2)^\alpha n^{d\alpha}}=1.$$
Hence we need to show that we have 
\eq{}\label{eq:ice}\lim_{\eps\downarrow 0}\lim_{{n\to+\infty}}\frac{1}{n^d}\sum_{x\in \T_n^d}\left| \sum_{z\in \Z^d\setminus \{0\}} \frac{\widehat \phi_\eps(z)\e^{-2\pi \imath z\cdot x} }{\|z\|^2} \widehat{f_n}(z)\right|^{\alpha}
= \int_{\T^d} \left| \sum_{z\in \Z^d\setminus \{0\}} \frac{\e^{-2\pi \imath z\cdot x}}{\|z\|^2} \widehat{f}(z)\right|^{\alpha} \De x.\eeq{}
We need this preliminary Lemma:
\begin{lemma}\label{lem:fourier}
There exists $C>0$ depending only on $f$ such that for all $n\in \N$
\[
\left|\widehat f(z)-\widehat{f_n}(z)\right|\le C n^{-1}.
\]
\begin{proof}
We can write 
\begin{align*}
\widehat f(z)-\widehat{f_n}(z)&= \sum_{x\in \T_n^d} \int_{B(x,\,\nicefrac{1}{2n})}\left[f( u)\cos\left( 2\pi z\cdot u\right)-f(x)\cos\left(2\pi  z\cdot x\right)\right]\De u\\
&+\imath\sum_{x\in \T_n^d} \int_{B(x,\,\nicefrac{1}{2n})}\left[f( u)\sin\left(- 2\pi z\cdot u\right)-f(x)\sin\left(-2\pi  z\cdot x\right)\right]\De u.
\end{align*}
Hence $\left|
\widehat f(z)-\widehat{f_n}(z)\right|$ is bounded above by the modulus of the two terms on the right-hand side of the previous equation. We will bound the first one, as the second is very similar. Using that the function $\psi:\,u\mapsto f(u)\cos\left(2\pi  z\cdot u\right)$ is $C^\infty(\T^d)$, we have from Taylor's series that
\begin{align*}
\left|f( u)\cos\left( 2\pi z\cdot u\right)-f(x)\cos\left(2\pi  z\cdot x\right)\right|\le \sup_{w\in \T^d}\left|\partial^\beta \psi(w)\right|\|x-u\|\le C n^{-1},
\end{align*}
where $\beta$ is a multi-index of degree $1$. Hence the conclusion follows.
\end{proof}
\end{lemma}
Let us now go back to \eqref{eq:ice}. Its left-hand side can be rewritten as
\[
\sum_{x\in \T_n^d}\left|\frac{1}{n^{\nicefrac{d}{\alpha}}} \sum_{z\in \Z^d\setminus \{0\}} \frac{\widehat \phi_\eps(z)\exp(-2\pi \imath z\cdot x) }{\|z\|^2} \widehat{f_n}(z)\right|^{\alpha}=:\sum_{x\in \T_n^d}\left|l_n^{(3)}(x)\right|^\alpha.
\]
As in the previous steps we write $l_n^3(x)=C_n^{(4)}(x)+l_n^{(4)}(x)$ with
\[
l_n^{(4)}(x):=\frac{1}{n^{\nicefrac{d}{\alpha}}} \sum_{z\in \Z^d\setminus \{0\}} \frac{\widehat \phi_\eps(z)\exp(-2\pi \imath z\cdot x) }{\|z\|^2} \widehat f(z)
\]
and
\[
C_n^{(4)}(x):=\frac{1}{n^{\nicefrac{d}{\alpha}}} \sum_{z\in \Z^d\setminus \{0\}} \frac{\widehat \phi_\eps(z)\exp(-2\pi \imath z\cdot x) }{\|z\|^2}\left( \widehat{f_n}(z)-\widehat f(z)\right).
\]
We need again to show that $n^{\nicefrac{d}{\alpha}}\left\|C_n^{(4)}\right\|_2$ goes to $0$. In order to do so, Lemma~\ref{lem:fourier} yields
\begin{align*}
\left\|C_n^{(4)}\right\|_2^2&={n^{-\nicefrac{2d}{\alpha}}}\sum_{z\in\Z_n^d\setminus\{0\}}\frac{\left|
\widehat{\phi_\eps}(z)
\right|^2\left| \widehat{f_n}(z)-\widehat f(z)\right|^2}{\|z\|^4}\le  \frac{C}{n^{\frac{2d}{\alpha}+2}}\sum_{z\in\Z_n^d\setminus\{0\}}\frac{\left|
\widehat{\phi_\eps}(z)
\right|^2}{\|z\|^4}\le  \frac{C}{n^{\frac{2d}{\alpha}+2}}.
\end{align*}
Here we have used the fast decay of $\phi_\eps$ at infinity.
Hence we get 
$$\lim_{n\to+\infty}\frac{1}{n^d}\sum_{x\in \T_n^d}\left| \sum_{z\in \Z^d\setminus \{0\}} \frac{\widehat \phi_\eps(z)\e^{-2\pi \imath z\cdot x} }{\|z\|^2} \widehat{f}(z)\right|^{\alpha}
= \int_{\T^d} \left| \sum_{z\in \Z^d\setminus \{0\}} \frac{\widehat \phi_\eps(z)\e^{-2\pi \imath z\cdot x}}{\|z\|^2} \widehat{f}(z)\right|^{\alpha} \De x.$$
Now noting that $f$ is a smooth function on $\T^d$ and $\left|\widehat f(z)\right|\le (1+\|z\|)^{-d-s}$ for $s\ge 0$ we can apply the dominated convergence theorem for $\eps \to 0$ and observing that $\widehat \phi_\eps(z)\to 1$ we obtain the result.
\end{proof}
\subsubsection{Scaling limit for regularly varying functions}
In this section we consider the scaling limit for a more general class of random variables. Since we are seeking a central limit type result it is natural to consider variables belonging to the domain of attraction of $\alpha$-stable distributions. 

Let $(\sigma(x))_{x\in \Z^d}$ be i.i.d. random variables satisfying Definition~\ref{firstassumption}; we can now start the proof of Theorem~\ref{thm:scaling_RV}.
\begin{proof}[Proof of Theorem~\ref{thm:scaling_RV}]An argument analogous to the one leading to \eqref{eq:for_later} shows that, by Proposition~\ref{prop:1.3LMPU} and the zero-mean property of test functions, it suffices to show that the statement holds for the field
\[
\la  \Xi_n, f \ra :=4\pi^2 n^{d-\frac{d}{\alpha}-2}   \sum_{z\in \T_n^d} w_\sigma(nz) \int_{B\left(z, \frac1{2n}\right)} f(t) \De t
\] 
where $$w_{\sigma}(z'):=(2d)^{-1} \sum_{x\in \Z_n^d} g(x,\,z') \sigma(x).$$ 
Let $(\rho(x))_{x\in \Z^d_n}$ be independent and distributed as $\rho_\alpha$ in \eqref{eq:doa:stable}. Then set
$$\la \widetilde \Xi_n, f \ra := 4\pi^2 n^{d-\frac{d}{\alpha}-2}   \sum_{z\in \T_n^d} w_{\rho}(nz) \int_{B\left(z, \frac1{2n}\right)} f(t) \De t$$
with $w_\rho$ defined as $w_\sigma$ replacing $\sigma$ by $\rho$. The proof will follow from Proposition~\ref{prop:slutsky} which will show the following equality in law:
$$\la \Xi_n, f \ra = \la \widetilde \Xi_n, f \ra + R_n$$
where $R_n$ goes to $0$ in probability. Thus it follows from Theorem~\ref{thm:scaling_RV_stable} that $\la \Xi_n, f\ra$ converges in distribution to $\la \widetilde \Xi_n, f \ra $ for all test functions $f$.
\end{proof}
To prove Proposition~\ref{prop:slutsky} we need to recall the following result. Consider a collection $(U_x)_{x\in \Z^d}$ of i.i.d. $U(0,1)$ random variables and let $(Y_x)_{x\in \Z^d}$ be a collection of i.i.d. random variables distributed as $\rho_\alpha$. We have 
\begin{lemma}[\protect{\citet{SimSto}, \citet[Lemma~3.3]{KluMik}}]\label{lemma:coupling}
Under the assumption of~\eqref{eq:doa:stable} we have that 
$$\left\{F_{\rho}^{\leftarrow}(U_x)\right\}_{x\in \Z^d}\overset{d}= (\rho(x))_{x\in \Z^d}\qquad \left\{ F_{\sigma(0)}^{\leftarrow}(U_x)\right\}_{x\in \Z^d}\overset{d}= (\sigma(x))_{x\in \Z^d}$$
and 
\begin{equation}
\lim_{n\to+\infty}{n^{-\frac{d}{\alpha}}} \sum_{x\in \Z_n^d}\left| \left[ F_{\rho}^{\leftarrow}(U_x)-F_{\sigma(0)}^{\leftarrow}(U_x)\right]\right|= 0
\end{equation}
in probability.
\end{lemma}
Now without of loss of generality we assume that $(\rho(x))_{x\in \Z^d}$ and $(\sigma(x))_{x\in \Z^d}$ live on the same probability space as in Lemma~\ref{lemma:coupling}. Let us now complete the proof of Theorem~\ref{thm:scaling_RV} by giving the proof of the last Proposition needed for it.

\begin{proposition}\label{prop:slutsky}
For $f\in C^\infty(\T^d)$ with mean zero, for every $\eps>0,$
$$\lim_{n\to+\infty}\prob\left( \left|\la \Xi_n, f \ra - \la \widetilde \Xi_n, f \ra\right|\ge \eps\right)= 0.$$
\end{proposition}
\begin{proof}
To obtain the above statement note that we have
$$\la \Xi_n, f \ra - \la \widetilde \Xi_n, f \ra= \sum_{x\in \Z_n^d} k_n(x) \left[ \sigma(x)-\rho(x)\right].$$
Here we have employed \eqref{eq:def_H}, \eqref{eq:def_k_n}.
Observe that 
\eq{}\label{eq:bound_k}\sup_{n\ge 1}\sup_{x\in \Z_n^d}| k_n(x)|\le Cn^{-d/\alpha}.\eeq{}
To prove this, we use an important technical estimate from \citet[Lemma~13]{CHR2016}: there exists $\mathcal M>0$ such that
$$\sum_{z\in \Z_n^d}\left|\widehat{L_n}(z)\right|\le \mathcal M n^{-d},$$
where $L_n$ is as defined in~\eqref{def:L_n}. Also $|\lambda_z|\ge C\|z\|^2 n^{-2}\ge Cn^{-2}$ for $\|z\|\ge 1$. Hence we get that
$$|k_n(x)| \le C n^{d-\frac{d}{\alpha}-2} \sum_{z\in \Z_n^d\setminus \{0\}} \left|\widehat{L_n}(z)\right| n^{2} \le C n^{-\frac{d}{\alpha}}.$$
Hence this proves \eqref{eq:bound_k}. 
Now we obtain Proposition~\ref{prop:slutsky} from Lemma~\ref{lemma:coupling}. 
\end{proof}

\appendix
\section{Appendix}\label{appendix}
\subsection{Proof of Lemma~\ref{lemma:5_points}}\label{sec:app:Lemma}
The proof of the Lemma requires a control on the tail behavior of the odometer series in the following way:
\begin{lemma}\label{claim:appendix}
Let $(Z_j)_{j\in \N}$ be $RV_{-\alpha}$, $\alpha\in (1,\,2)$. Moreover $\E[Z_j]=0$ holds for all $j$ and $\prob(Z_j>t)=\prob(Z_j<-t)$ for all $t\ge 0$. Let $(c_j)_{j\in \N}$ satisfy
\[
\sum_{j\ge 1}\left|c_j\right|^\delta<+\infty,\;\text{for some }\delta<\alpha.
\]
Then for any $M\ge 1$ there exist $n_1=n_1(M)$, $a>0$ for which
\eq{}\label{eq:sum_cond}
\prob\Bigg( \Big|\sum_{j\ge n_1}c_j Z_j\Big|>\frac1M\Bigg)\le M^{-a}.
\eeq{}
\end{lemma}
\begin{proof}
Let $\delta$ be as in the assumptions. Since $Z_j$ is $RV_{-\alpha}$ there exists $x_1$ such that for all $x\ge x_1$ we have $\prob(|Z_1|>x)\le \nicefrac{1}{2}\, x^{-\delta}$.
We use Karamata's theorem \cite[Theorem~0.6]{Resnick} which gives that
$$U(x):= \E\left[ |Z_1|^2\one_{|Z_1|\le x}\right]\in RV_{2-\alpha}.$$
Hence there exists $x_2$ such that 
\eq{}\label{eq:U}
U(x)\le \frac{1}{2}\, x^{\alpha-\delta},
\eeq{}
for all $x\ge x_2$. 

Fix $\eps>0$. The following conditions hold for $n_1$ large enough:
\begin{enumerate}[label=(C.\arabic*),ref=(C.\arabic*)]
\item \label{item:C.1} $\sum_{j=n_1}^{+\infty} |c_j|^{\delta} < \eps^{2{\delta}}$, 
\item \label{item:C.2}$\eps |c_j|^{-1} \ge \max\{x_1, x_2\}$ and $|c_j|\le 1$, $j>n_1$. 
\end{enumerate}
Note that such choices can be made as $c_j\to 0$ as $j\to +\infty$. We have then
\begin{align*}
\prob\Bigg( \Big|\sum_{j=n_1}^{+\infty} c_j Z_j\Big|>\eps\Bigg)&\le \prob\Bigg( \Big|\sum_{j=n_1}^{+\infty} c_j Z_j\Big|>\eps, \sup_{j\ge n_1} |c_j Z_j|>\eps\Bigg)\\
&+\prob\Bigg( \Big|\sum_{j=n_1}^{+\infty} c_j Z_j\Big|>\eps, \sup_{j\ge n_1} |c_j Z_j|\le \eps\Bigg)\\
&\le \sum_{j=n_1}^{+\infty} \prob( |c_j Z_j|>\eps) + \prob\Bigg( \Big|\sum_{j=n_1}^{+\infty} c_j Z_j\one_{\left\{|c_j Z_j|\le \eps\right\}}\Big|>\eps\Bigg).
\end{align*}
First we tackle the first sum. Note that
$$\sum_{j\ge n_1} \prob( |c_j Z_j|>\eps) = \sum_{j\ge n_1} \prob\left( |Z_j|>\frac{\eps}{|c_j|}\right)\le \frac{\eps^{-\delta}}{2} \sum_{j\ge n_1} |c_j|^{\delta}< \frac{\eps^{\delta}}{2}$$
thanks to \ref{item:C.1}.
Next we handle the second term with Markov's inequality:
$$P\Bigg( \Big|\sum_{j=n_1}^{+\infty} c_j Z_j\one_{\left\{|c_j Z_j|\le \eps\right\}}\Big|>\eps\Bigg)\le \eps^{-2}\E\Bigg[\Big| \sum_{j=n_1}^{+\infty} c_j Z_j \one_{\left\{|Z_j|\le \frac{\eps}{|c_j|}\right\}}\Big|^2\Bigg].$$
Let us denote by $W_j:= c_j Z_j \one_{\left\{|Z_j|\le \eps/|c_j|\right\}}$. Now note that the independence of the $Z_j$'s, Fatou's lemma and the monotone convergence theorem imply
$$\E\Bigg[\Big|\sum_{j=n_1}^{+\infty} W_j\Big|^2\Bigg]\le \sum_{j=n_1}^{+\infty} \E\left[ W_j^2\right] +\left( \sum_{j=n_1}^{+\infty} \E[|W_j|]\right)^2.$$
We bound each one of the terms above.  Observe that $\E\left[ W_j^2\right] =c_j^2 U( \eps/|c_j|)$. Since $\eps|c_j|^{-1}>x_2$ by \ref{item:C.2} we have that from \eqref{eq:U} 
$$\E\left[ W_j^2\right]\le \frac12\eps^{\alpha-\delta}|c_j|^{\alpha-\delta}$$
follows
Hence we have
$$\sum_{j=n_1}^{+\infty} \E\left[ W_j^2\right] \le \frac{\eps^{2-\delta}}{2}\sum_{j=n_1}^{+\infty} |c_j|^{\delta} \le \frac12\eps^{2+\delta}.$$
Now an argument analog to \citet[Equation~(2.6)]{KokTaq} gives us
$$\E\left[ |W_j|\right] \le \frac{(1+\alpha-\delta)\delta}{\delta-1}\eps^{1-\delta} |c_j|^{\delta}.$$
So we get that for some constant $C>0$
$$\sum_{j=n_1}^{+\infty} \E\left[|W_j|\right] \le \frac{(1+\alpha-\delta)\delta}{\delta-1}\eps^{1-\delta} |c_j|^{\delta}\sum_{j=n_1}^{+\infty} |c_j|^{\delta} \le C \eps^{1+\delta}.$$
This shows that
$$\prob\left( \left|\sum_{j=n_1}^{+\infty} c_j Z_j\right|>\eps\right) \le \max\left\{ \eps^{\delta}, \frac12\eps^{2+\delta}, C \eps^{1+\delta}\right\}.$$
This completes the proof.
\end{proof}
\begin{corollary}\label{cor:truncated_series} For all $M\ge 1$ there exist $n_1$ and $N\ge n_1$ such that
$$\prob\left(\left| \sum_{j=n_1}^{N} c_j Z_j\right|\ge \frac1M\right)\le M^{-a}.$$
\end{corollary}
\begin{proof}
It follows from Lemma~\ref{claim:appendix} setting $c_j:=0$ for all $j>N_1$.
\end{proof}
\begin{proof}[Proof of Lemma~\ref{lemma:5_points}] \noindent
\begin{enumerate}[label=(\Roman*),leftmargin=*,rightmargin=-.01pt]
\item The series is finite almost surely by \citet[Theorem~2.1 ii)]{Cline} and \ref{item:(b)}.
\item The proof follows the steps of \citet[Lemma~5.5 b), d)]{LMPU}. While d) carries over to our setting, we have a substantial difference in b), where we do not have finite variance of the random variables
\[
v_{\gamma,\,N}:=\sum_{j=1}^N g(o,\,\gamma y_j)Y_{\gamma y_i},
\] 
for $\gamma\in \Gamma$.
However, we can estimate $\prob(|v_{\gamma,\,N}-v_{e,N}|>\eps)$, $N\in \N$, $\eps>0$ by Corollary~\ref{cor:truncated_series} and obtain the same conclusion.
\item Choose $\eps_1\in (1,\,\alpha)$. Since $L$ is slowly varying, we have that $t^{\eps_1}L(t)\to +\infty$ as $t\to+\infty$. Hence there exists a $t_0$ such that $L(t)>t^{-\eps_1}$ for $t\ge t_0$, and so 
\eq{}\label{eq:t_large}
P(Y_o<-t)> t^{-(\alpha+\eps_1)}>0,\quad t\ge t_0.
\eeq{}
Choose $M\ge 1$ arbitrarily large. 
We use Lemma~\ref{claim:appendix} for $c_j:=r^{-1}g(o,\,y_j)$ and $Z_j:=Y_{y_j}$ to find an $n_1=n_1(M)$ such that
\eq{}\label{eq:tail_M_sum}
\prob\Bigg(\frac{1}{r}\sum_{i\ge n_1}g(o,\,y_i)Y_{y_i}> M\Bigg)\le M^{-a}.
\eeq{}
Observe furthermore that on the event $\left\{ Y_{y_i}<-t:\,i\le n_1-1\right\}$ one has
\eq{}\label{eq:part_1}
\frac{1}{r}\sum_{i\le n_1-1}g(o,\,y_i)Y_{y_i}\le -\frac{t}{r}\sum_{i\le n_1-1}g(o,\,y_i).
\eeq{}
Moreover we can choose $t=t(M)\ge t_0$ large enough so that
\eq{}\label{eq:part_2}
\frac{t}{r}\sum_{i\le n_1-1}g(o,\,y_i)>2M.
\eeq{}
Thus for $t=t(M)$, $n_1=n_1(M)$ as above
\begin{align*}
\prob\left(v_e(o)<-M\right)&\ge\prob\Bigg(\frac{1}{r}\sum_{i\le n_1-1}g(o,\,y_i)Y_{y_i}<-2M\Bigg)\prob\Bigg(\frac{1}{r}\sum_{i\ge n_1}g(o,\,y_i)Y_{y_i}<M\Bigg)\\
&\stackrel{\eqref{eq:part_1},\eqref{eq:part_2}}{\ge}\prob\Bigg(Y_{y_i}<-t:\,i\le n_1-1\Bigg)\prob\Bigg(\frac{1}{r}\sum_{i\ge n_1}g(o,\,y_i)Y_{y_i}\le M\Bigg)\\
&\stackrel{\eqref{eq:t_large},\,\eqref{eq:tail_M_sum}}{\ge} t^{-(\alpha+\eps_1)N}\left(1-M^{-a}\right)>0.
\end{align*}
Hence by ergodicity of $v_e$ and the fact that $M$ is arbitrary, we have that 
$$\prob\left( \inf_{x\in V} v_e(x)<-t\right)=1.$$  
\end{enumerate}
\end{proof}
\subsection{Stable distributions}\label{subsec:app:stable}
We have shown that the characteristic functional of $\la\Xi_n,\,f\ra$ has the form $\exp(-\mathcal L_\alpha(f))$, where 
$$\mathcal L_\alpha(f)= \int_{\T^d} \left| \sum_{z\in \Z^d\setminus \{0\}} \frac{\exp(-2\pi \imath z\cdot x)}{\|z\|^2} \widehat{f}(z)\right|^{\alpha} \De x.$$ 
We want to investigate properly the measure associated to the latter characteristic functional. Recall the definition of the space $\mathcal T:=C^\infty(\T^d)/{\sim}$. This is a nuclear space and it is reflexive (by \citet[Section~8.4.7]{EdwFun} and the fact that the quotient of a reflexive space by
a closed subspace is reflexive). We would like here to show that this functional defines a measure on $\mathcal T^*$ via the Bochner-Minlos theorem. If this is true, then
\[
(-\Delta)^{-1}:\,\mathcal T^{**}=\mathcal T\to L_\alpha(\T^d)
\]
defines an $\alpha$-stable measure on $\mathcal T^{*}$ (cf. \citet[Theorem~5]{Linde} in the setting of Banach spaces). 
\begin{theorem}[Bochner-Milnos]
Let $V$ be a nuclear space. Then a complex valued function $\Phi$ on $V$ is the characteristic function of a probability measure $\nu$ on $V^*$ if and only if $\Phi(0)=1$, $\Phi$ is continuous and $\Phi$ is positive definite, that is,
$$\sum_{j,\,k=1}^n z_j \overline{z_k} \Phi(v_j-v_k)\ge 0$$
for all $v_1,\,\ldots,\, v_n\in V$ and $z_1, \ldots, z_n\in \mathbb C$.
\end{theorem}
We apply Bochner-Minlos theorem to obtain
\begin{theorem}
The functional $\Phi(f):=\exp\left(-\mathcal L_\alpha(f)\right)$ on the space $\mathcal T$ is the characteristic function of a probability measure on $\mathcal T^*$.
\end{theorem}
\begin{proof}
From Bochner-Minlos theorem we need to check three assumptions.
\begin{enumerate}[leftmargin=*,rightmargin=-.01pt]
\item Recall
$$(-\Delta)^{-1}f(x)= \sum_{z\in \Z^d\setminus\{0\}} \frac{\widehat f(z)}{\|z\|^{2} }\exp(-2\pi \imath z\cdot x).$$
Using $\left|\e^{-x}-\e^{-y}\right|\le |x-y|$ we obtain for two arbitrary $f_1,\,f_2\in C^\infty(\T^d)$
\begin{align*}
\left|\e^{-\mathcal L_\alpha(f_1)}-\e^{-\mathcal L_\alpha(f_2)}\right|&\le \left|\mathcal L_\alpha(f_1)-\mathcal L_\alpha(f_2)\right|\\
&=\left|\int_{\T^d}\left|(-\Delta)^{-1}f_1(x)\right|^\alpha \De x-\int_{\T^d}\left|(-\Delta)^{-1}f_2(x)\right|^\alpha \De x\right|.
\end{align*}
From \citet[Lemma 4.7.2]{ST} we see that the last term is bounded above by
\begin{align*}
2^{\nicefrac{1}{\alpha}}\alpha\left(\|(-\Delta)^{-1}f_1\|_\alpha^{\alpha-1}+\|(-\Delta)^{-1}f_2\|_\alpha^{\alpha-1}\right)&\left(\int_{\T^d}\left|(-\Delta)^{-1}(f_1-f_2)(x)\right|^\alpha\De x\right)^{\nicefrac{1}{\alpha}}\\
&=:C_\alpha \left\|(-\Delta)^{-1}(f_1-f_2)\right\|_{\alpha}.
\end{align*}
One case see that $(-\Delta)^{-1}f\in L^\alpha(\T^d)$ if $f$ is smooth: in fact
\[
\left|(-\Delta)^{-1}f(x)\right|\le \sum_{z\in \Z^d\setminus\{0\}} \frac{\widehat{f}(z)}{\|z\|^2} <+\infty
\]
due to the fact that $\widehat f(0)=0, \,\|z\|\ge 1$ and by the decay properties of $\widehat f$ \cite[Theorem 5.4]{Roe}. Then $(-\Delta)^{-1}f\in L^\infty(\T^d)$ and so is in any $L^\alpha$. So we notice now that
\begin{align*}
&\left\|(-\Delta)^{-1}(f_1-f_2)\right\|_{\alpha}\le  \left\|(-\Delta)^{-1}(f_1-f_2)\right\|_{2}\\
&\le
 \Bigg(\sum_{z\in \Z^d\setminus\{0\}}{\|z\|^{-4}}{\left|\widehat{f_1}(z)-\widehat{f_2}(z)\right|^2}\Bigg)^{\nicefrac{1}{2}}\le\Bigg(\sum_{z\in \Z^d}{\left|\widehat{f_1}(z)-\widehat{f_2}(z)\right|^2}\Bigg)^{\nicefrac{1}{2}}
\end{align*}
using the orthogonality of the characters in the second-to-last equality and the fact that $\|z\|>1$ in the last. Parseval's theorem yields then
\[
\left\|(-\Delta)^{-1}(f_1-f_2)\right\|_{\alpha}\le \left(\int_{\T^d}\left({f_1}(x)-{f_2}(x)\right)^2\De x\right)^{\nicefrac{1}{2}}\le \sup_{x\in \T^d}\left|{f_1}(x)-{f_2}(x)\right|.
\]
Since the Fr{\'e}chet topology on $C^\infty$ is given by the uniform convergence of all derivatives, we have continuity.
\item The fact that $\Phi(0)=1$ is immediate.
\item The positive definiteness of $\exp( -\mathcal L_\alpha(f))$ follows since it is a limit of positive definite functionals.
\end{enumerate}
\end{proof}
\bibliographystyle{abbrvnat}
\bibliography{literaturASP}
\end{document}